\documentclass[11pt,reqno]{amsart}
\usepackage{hyperref}
\usepackage{url}
\usepackage{amssymb}
\usepackage[T1]{fontenc}
\usepackage{extarrows}
\usepackage{lmodern}
\usepackage[utf8]{inputenc}
\usepackage{ifthen}
\usepackage{mathtools}
\usepackage{xcolor}
\usepackage{amsmath}
\usepackage{booktabs}
\usepackage{array}
\usepackage{geometry}
\usepackage{tikz}
\newtheorem{thm}{Theorem}
\newtheorem{coro}{Corollary}

\newtheorem{lem}{Lemma}
\newtheorem{rem}{Remark}
\newtheorem{prop}{Proposition}

\theoremstyle{definition}
\newtheorem{definition}{Definition}
\theoremstyle{remark}

\numberwithin{equation}{section}

\newcommand{\R}{\mathbb{R}}

\newcommand{\I}{\mathcal{I}}

\usepackage{tikz,enumerate}
\begin{document}
\title{Codes and designs in multivariate $Q$-polynomial association schemes}
\author{Minjia Shi}
\author{Jing Wang}
\author{ Patrick Sol\'e}
\footnotetext[1]{This research is supported by the National Natural
Science Foundation of China under Grants W2611002 and 12471490.
(\textit{Corresponding author: Minjia Shi})}
\footnotetext[2]{Minjia Shi is with the Key Laboratory of Intelligent
Computing and Signal Processing, Ministry of Education,
School of Mathematical Sciences, Anhui University, Hefei, China
(email: smjwcl.good@163.com).}
\footnotetext[3]{Jing Wang is with the School of Mathematical Sciences,
Anhui University, Hefei 230601, China
(email: wangjing031010@163.com).}
\footnotetext[4]{Patrick Solé is with I2M
(Aix Marseille Univ, CNRS), Marseille, France
(email: sole@enst.fr).}
\maketitle
\begin{abstract}
We generalize the fundamental bounds of Delsarte's thesis (1973) on codes of given degree and designs of given strength in the new setting of Bannai et al. (2025). We assume the scheme is weakly metric in the sense of (Sol\'e, 1989). We give upper bounds on the size of codes of given degree, and also on the size of codes with a given number of pairwise distances. Codes meeting these bounds are characterized by the relation of suitable annihilators with the degree (resp. distance) Wilson polynomial. We give two analogues of the Rao bound on the size of designs with given strength. Designs meeting that bound we call degree tight designs or distance tight designs depending on the bound met.
Applications to the Lee distance, mixed level orthogonal arrays, ordered orthogonal arrays, and more are given. The formal duality between codes and designs, connecting perfect codes and tight designs, is made concrete in
self-dual translation schemes.
\end{abstract}


\section{Introduction}
The landmark PhD thesis of Delsarte in 1973 created the field of algebraic combinatorics  \cite{D2}. In that seminal work, the idea of {\it association schemes} (an algebraic structure coming from both  statistics \cite{Ba} and group theory \cite{BI,BCN}) and {\em linear programming} (a tool from combinatorial optimization \cite{Ch}) are combined
to derive bounds on the size of certain combinatorial objects (codes, designs), and characterize the special objects meeting these bounds (perfect codes, tight designs). The applicability of these techniques was limited to some particular association schemes
(metric schemes, Q-polynomial association schemes) which enjoy a very simple generator structure of their Bose-Mesner algebra. However, in many applications, association schemes that are not metric and/or not Q-polynomial are needed \cite{AGKP,NR,MW,NRT,qary,S3}. In \cite{S2} was introduced the notion of {\em weakly metric} association scheme, which was further developed in \cite{S1}. In 2025, the concept of {\em multivariate} P-polynomial (and also Q-polynomial) association schemes was introduced by Bannai et al. \cite{B+}, generalizing \cite{BC}. The related notion of an \(m\)-distance-regular graph, generalizing metric schemes and distance-regular graphs, was developed by Bernard et al. \cite{B++}.

This paper is an attempt to lay the foundations of {\em multivariate} Delsarte theory of codes and designs in {\em general} association schemes. It is clear that any association scheme is $\ell$-variate Q-polynomial for some $\ell.$ See \cite[Remark 3, p.10]{B+} for the P-polynomial version of this observation.
Thus, we consider extremal properties of codes and designs in multivariate Q-polynomial association schemes. A result of \cite{D2} was popularized in \cite{D}. It is an upper bound on the size of a set with a given number of pairwise relations, the so called degree of the set. In \cite{D}, this bound is given for the Hamming scheme as an upper bound on the size of codes with a given number of pairwise distances. In the situation of a weakly metric Q-polynomial association
scheme there are two possible generalizations of this bound: for codes with a given degree or for codes with a given number of distances. This happens because the scheme being not metric in general, several relations correspond to the same distance. We give these two possible generalizations. The pairwise relations of codes in the first case are
zeros of the degree Wilson polynomial. The pairwise distances of codes meeting the second bound are zeros of the distance Wilson polynomial.
Furthermore, in the special case of self-dual translation schemes, additive codes $Y$  attaining the relation-degree bound \(|Y|=M(r)\) have perfect duals. Moreover, the covering radius of such a perfect dual is restricted by the regularity of the dispersion function of the scheme.

In an interlude section we show that many multivariate Q-polynomial association schemes of interest in Combinatorics and Information Theory satisfy a hypothesis needed for the second bound. Thus we consider bivariate (q-ary Johnson scheme \cite{qary}), multivariate schemes like the Lee scheme \cite{S3}, mixed alphabet scheme (related to mixed level orthogonal arrays \cite{DP}) and NRT scheme \cite{NR,NRT}.

In the final part of the paper, we consider designs in a weakly metric multivariate  Q-polynomial scheme.
The notion of $T$-design with $T$ a set makes sense in any multivariate Q-polynomial scheme. Note that designs in this general situation were explored in \cite{DP,qary}. To define $t$-designs with $t$ an integer we need a weakly metric structure. We derive a generalized Rao lower bound on the size of such designs. Designs meeting that bound with equality we call {\em distance} tight designs. To be complete we also give a lower bound on the size of $T$-designs. Tight designs in that context we call {\em degree} tight designs. If furthermore, the scheme is a self-dual  translation scheme then the dual of a distance tight design is shown to be a perfect code.

While all these results are natural generalizations of one-variable results, the extension is not immediate. Passing from one variable polynomials to multivariate polynomials requires technical assumptions like the degree filtration property attached to the spectral node set (\S 2.2).

The material is arranged as follows. Section 2 contains preliminary notions on association schemes needed for the rest of the paper. Section 3 derives an upper bound on codes of given degree. Section 4 derives a similar bound
on codes with a given number of pairwise distances. Section 5 gives examples of schemes where the previous bounds apply. Section 6 considers designs: Rao bound, tight designs and the Wilson polynomial.
Section 7 concludes the article. An appendix collects facts on and examples of self-dual translation association schemes.
\section{Preliminary}
\subsection{Weakly metric association schemes}

We recall the axioms of commutative association schemes. See \cite{BI,BCN, MS} for background and details.
\begin{definition}
Let $X$ be a finite set with at least two elements and let $R = \{R_0, R_1, \ldots, R_s\}$ be a family of $s + 1$ relations $R_i$ on $X$ for $s \geq 1$.
The pair $(X, R)$ is called an {\bf association scheme}  with $s$ classes if the following conditions are satisfied:
\begin{enumerate}
\item[(A1)] The set $R$ is a partition of $X^2$ and $R_0 = \{(x, x) \mid x \in X\}$ is the diagonal relation;
\item[(A2)] For $i = 0, 1, \ldots, s$, the inverse $R_i^{-1} = \{(y, x) \mid (x, y) \in R_i\} $ of the relation $R_i$ also belongs to $R$;
\item[(A3)] For any triple of integers $i, j, k = 0, 1, \ldots, s$, there exists a number $ p_{ij}^k = p_{ji}^k$ such that
$$ |\{z \in X \mid (x, z) \in R_i, (z, y) \in R_j\}| = p_{ij}^k $$
for all \( (x, y) \in R_k \).
\end{enumerate}
\end{definition}

An association scheme is \textbf{symmetric} if \(R_i^{-1}=R_i, 0\leq i\leq s\). Equivalently, each adjacency matrix is real symmetric. Throughout this paper, unless otherwise stated, all association schemes are assumed to be symmetric.

We say that an association scheme $(X,R)$ is {\bf weakly metric} (shortly WMAS) for a quasi-distance $d$ if the quasi-distance is constant on the classes of the scheme \cite{S2}. 
More precisely, let \(\mathcal{I}\) be the index set of relations, with \(0\) indexing the diagonal relation. A quasi-distance \(d:X\times X\to\mathbb R_{\geq 0}\), satisfying \(d(x,x)=0\) and the triangle inequality, makes the scheme weakly
metric if there exists a function \(d_\mathcal{I}:\mathcal{I}\to\mathbb R_{\geq 0}\) such that
\[
d(x,y)=d_\mathcal{I}(\alpha),~\text{whenever }(x,y)\in R_\alpha.
\]
Throughout this paper, we assume that \(d_\mathcal{I}(0)=0\) and \(d_\mathcal{I}(\alpha)>0\) for every \(\alpha\ne 0\).
We use \(d_\mathcal{I}\) for the relation index distance and \(d_{\mathcal D^*}\) for the weight on primitive idempotent indices, specified in the relevant sections.

For each relation $R_i$ ($i = 0, 1, \dots, s$), we define its adjacency matrix $A_i$ by:
\[
A_i(x, y) =
\begin{cases}
1 & \text{if } (x, y) \in R_i, \\
0 & \text{if } (x, y) \notin R_i.
\end{cases}
\]
The set of all complex linear combinations of these adjacency matrices forms a complex algebra:
\[
\mathfrak{A} = \left\{ \sum_{i=0}^s \alpha_i A_i \,\bigg|\, \alpha_i \in \mathbb{C} \right\},
\]
known as the {\bf Bose-Mesner algebra} \cite{BI,BCN,MS} of the association scheme $(X,R)$. The Bose-Mesner algebra admits two natural bases. The first is the set of adjacency matrices $\{A_i\}_{i=0}^s$. The second is a unique basis of primitive idempotents $\{E_j\}_{j=0}^s$ satisfying:
\[
E_j^2 = E_j,~E_j E_k = 0 \ (j \neq k),~\sum_{j=0}^s E_j = I_{|X|},
\]
where $I_{|X|}$ denotes the identity matrix of order $|X|$. The linear transformation between these bases is given by:
\[
A_i = \sum_{j=0}^s P_i(j) E_j,~i = 0, 1, \dots, s.
\]
The complex numbers $P_i(0), P_i(1), \dots, P_i(s)$ are the eigenvalues of $A_i$, also called the {\bf first eigenvalues} of the association scheme.
Similarly, the complex numbers \( Q_j(0), Q_j(1), \ldots, Q_j(s) \) where the second eigenmatrix
\(Q = (Q_j(i))_{i,j}\) is defined by
\[
|X|E_j = \sum_{i=0}^s Q_j(i)A_i,~j = 0, 1, \ldots, s,
\]
are called the {\bf second eigenvalues} of the association scheme.
Let \(v_i\) denote the valency of the relation $R_i$, and let \(\mu_j=\operatorname{rank}(E_j)\) denote the multiplicity of the primitive idempotent $E_j$. For the symmetric association schemes considered below, we shall use the
following standard identities \cite[Theorem~11.2.1]{uvp}.

{\lem\label{eigenmatrix-identities}
For $0\leq i,k,j,l\leq s$, the first and second eigenvalues satisfy
\begin{enumerate}
\item[(1)]
\[
P_0(j)=1,~Q_0(i)=1,~P_i(0)=v_i,~Q_j(0)=\mu_j.
\]
\item[(2)]
\[
P_i(j)P_k(j)=\sum_{l=0}^{s}p_{ik}^{\,l}P_l(j).
\]
\item[(3)]
\[
\mu_jP_i(j)=v_iQ_j(i),
\]
and
\[
\sum_{j=0}^{s}\mu_jP_i(j)P_k(j)=|X|v_i\delta_{ik},~\sum_{i=0}^{s}v_iQ_j(i)Q_l(i)=|X|\mu_j\delta_{jl}.
\]
\item[(4)]
\[
|P_i(j)|\leq v_i,~|Q_j(i)|\leq\mu_j.
\]
\end{enumerate}
Moreover,
\[
PQ=QP=|X|I.
\]
}

The Bose-Mesner algebra $\mathfrak{A}$ is also closed under entrywise multiplication, denoted by $\circ$ and called the {\bf Hadamard product} \cite{B+}.

\begin{definition}
Let $X$ be a finite additive abelian group, and let $\mathfrak{X} = \bigl(X, \{R_i\}_{i=0}^s\bigr)$ be a $s$-class association scheme, then $\mathfrak{X}$ is called a {\bf translation association scheme} \cite{TAS} if
\[
(x, y) \in R_i \iff (x+z, y+z) \in R_i,~\text{for all } z \in X \text{ and } i.
\]
\end{definition}
Let \(X_i = \bigl\{x \in X \mid (0, x) \in R_i\bigr\},~\text{for}~0 \leq i \leq s\), then $X_0, X_1, \dots, X_s$ give a partition of $X$, and \((x, y) \in R_i \iff y - x \in X_i~,~\text{for}~0 \leq i \leq s\).

\begin{definition}
An {\bf e-perfect code} $C \subset X$ in the WMAS $(X,R)$ is then defined as a code satisfying the partition
\[
X=\coprod_{c \in C} B(c,e),
\]
where $\coprod$ means disjoint union, and \(B(c,e)=\{ y \in X \mid d(y,c) \le e\}\) is the ball of radius e for the quasi-distance $d$ centered in $c.$
\end{definition}
Equivalently, if \(e\) is the packing radius of \(C\), then \(C\) is perfect if and only if
\[
|C|\cdot |B(c,e)|=|X|.
\]

\subsection{Multivariate Q-polynomial schemes}
We begin with general definitions on polynomials in several variables.
\begin{definition}
A \textbf{monomial order} \cite{B+} $\leq$ on $\mathbb{C}[x_1, x_2, \ldots, x_\ell]$ is a relation on the set of monomials $x_1^{n_1} x_2^{n_2} \ldots x_\ell^{n_\ell}$ satisfying:
\begin{enumerate}
\item[(1)] $\leq$ is a total order;
\item[(2)] for monomials $u, v, w$, if $u \leq v$, then $wu \leq wv$;
\item[(3)] $\leq$ is a well-ordering, i.e., any non-empty subset of the set of monomials has a minimum element under $\leq$.
\end{enumerate}
\end{definition}
Since each monomial $x_1^{n_1} x_2^{n_2} \ldots x_\ell^{n_\ell}$ can be associated to the tuple
$(n_1, n_2, \ldots, n_\ell) \in \mathbb{N}^\ell$, we use the same order and notation $\leq$ on
\[
\mathbb{N}^\ell = \{(n_1, n_2, \ldots, n_\ell) \mid n_i \text{ are non-negative integers}\}.
\]
Note that for $\ell = 1$, there is a unique monomial order, which is the one associated to the natural ordering of the integers.

\begin{definition}
For a monomial order $\leq$ on $\mathbb{N}^\ell$, an $\ell$-variate polynomial $v(x_1, x_2, \ldots, x_\ell) = v(x)$ is said to be of \textbf{multi-degree} \cite{B+} $n \in \mathbb{N}^\ell$ if $v(x)$ is of the form
\[
v(x) = \sum_{a \leq n} f_a x^a, \quad x^a = x_1^{a_1} x_2^{a_2} \ldots x_\ell^{a_\ell},
\]
with $f_a \in \mathbb{C}$ and $f_n \neq 0$.
\end{definition}

We shall say that a monomial order \(\leq\) on \(\mathbb{N}^{\ell}\) is
\(L^1\)-compatible,
\[
a\leq b~\Longrightarrow~|a|\leq |b|,~\text{for all}~a,b\in\mathbb{N}^{\ell},
\]
where \(|a|=a_1+\cdots+a_\ell\). Under such an order, a polynomial of multi-degree \(\gamma\) has ordinary total degree \(|\gamma|\).

\begin{definition}
Let $\mathcal{D}^* \subset \mathbb{N}^\ell$ containing $\epsilon_1, \epsilon_2, \dots, \epsilon_\ell$, and let $\leq$ be a {monomial order} on $\mathbb{N}^\ell$. A commutative association scheme $\mathfrak{X} = (X, R)$ with the primitive idempotents $\{E_j\}_{j \in J}$ is called {\bf$\ell$-variate $Q$-polynomial} \cite{B+} on the domain $\mathcal{D}^*$ with respect to $\leq$ if the following three conditions are satisfied:
\begin{enumerate}
\item[(1)] if $(n_1, \dots, n_\ell) \in \mathcal{D}^*$ and $0 \leq m_i \leq n_i$ for $i = 1, \dots, \ell$, then $(m_1, \dots, m_\ell) \in \mathcal{D}^*$;
\item[(2)] there exists a relabeling of the primitive idempotents:
\[
\{E_j\}_{j \in J} = \{E_\alpha\}_{\alpha \in \mathcal{D}^*},
\]
such that, for $\alpha \in \mathcal{D}^*$,
\[
|X| E_\alpha = v_\alpha^*\bigl( |X| E_{\epsilon_1}, |X| E_{\epsilon_2}, \dots, |X| E_{\epsilon_\ell} \bigr) (\text{under the Hadamard product}),
\]
where $v_\alpha^*(\mathbf{x})$ is an $\ell$-variate polynomial of multidegree $\alpha$ with respect to $\leq$, and all monomials $\mathbf{x}^\beta$ in $v_\alpha^*(\mathbf{x})$ satisfy $\beta \in \mathcal{D}^*$;
\item[(3)] for $i = 1, 2, \dots, \ell$ and $\alpha = (n_1, n_2, \dots, n_\ell) \in \mathcal{D}^*$, the product $E_{\epsilon_i} \circ E_{\epsilon_1}^{\circ n_1} \circ E_{\epsilon_2}^{\circ n_2} \circ \dots \circ E_{\epsilon_\ell}^{\circ n_\ell}$ is a linear combination of
\[
\left\{ E_{\epsilon_1}^{\circ m_1} \circ E_{\epsilon_2}^{\circ m_2} \circ \dots \circ E_{\epsilon_\ell}^{\circ m_\ell} \mid \beta = (m_1, m_2, \dots, m_\ell) \in \mathcal{D}^*, \beta \leq \alpha + \epsilon_i \right\}.
\]
\end{enumerate}
\end{definition}

Let $\mathcal{I}$ be the index set of relations. For each $\alpha\in\mathcal{I}$, define the corresponding {\bf spectral node} by \(z_\alpha=\bigl(Q_{\epsilon_1}(\alpha),\ldots,Q_{\epsilon_\ell}(\alpha)\bigr)\).
For $\gamma\in\mathcal{D^*}$, write $\Phi_\gamma=v_\gamma^*$, then
\(Q_\gamma(\alpha)=\Phi_\gamma(z_\alpha)\). The set \(Z=\{z_\alpha:\alpha\in \mathcal{I}\}\) is called the {\bf spectral node set } associated with the chosen multivariate $Q$-polynomial structure.

\textbf{Hypothesis $F_{\deg}$.}
We assume that the spectral node set $Z$ and the functions $\{\Phi_\gamma|_Z:\gamma\in\mathcal D^*\}$ satisfy the following property: for every $s\geq 0$, if $F(x)$ is a polynomial of ordinary total degree at most $s$, then its restriction to $Z$ satisfies
\[
F|_Z\in \operatorname{span}\{\Phi_\gamma|_Z:\gamma\in \mathcal{D^*},\ |\gamma|\leq s\}.
\]

\subsection{Linear programming in association schemes}
Let \(N = \{0,1,\ldots,n\}\) and let \(A = [A_k(i)]\) be a matrix of \(R(N,N)\) such that \(A_0(i)=1\) for all \(i\in N\) and \(A_k(0)>0\) for all \(k\in N\).
Let \(\mathcal{M}\subseteq N\) be a subset containing \(0\), denote \(\mathcal{M}^* = \mathcal{M}\setminus\{0\}\) and \(N^* = N\setminus\{0\}\).
The \textbf{linear-programming problem} \((A,\mathcal{M})\) \cite{D2} is defined as follows:
\[
(A,\mathcal{M})\qquad
\begin{cases}
\text{max}\; g = \displaystyle\sum_{i\in \mathcal{M}} b_i, \\[1em]
\displaystyle\sum_{i\in \mathcal{M}} b_i A_k(i) \;\ge\; 0, & \forall k\in N^*, \\[1em]
b_i \;\ge\; 0, & \forall i\in \mathcal{M}^*.
\end{cases}
\tag{2.1}
\]

An \(n+1\)-tuple \(b = (b_0,b_1,\ldots,b_n)\) is called a \textbf{program} of \((A,\mathcal{M})\) if it satisfies \(b_0=1\), \(b_i=0\) for \(i \in N-\mathcal{M}\), and the constraints (2.1). A \textbf{maximal program} is a program that maximizes \(g\); its value is denoted by \(g(A,\mathcal{M})\).

The \textbf{dual problem} \((A,\mathcal{M})'\) is given by
\[
(A,\mathcal{M})'\qquad
\begin{cases}
\text{min}\; \gamma = \displaystyle\sum_{k\in N} \beta_k A_k(0), \\[1em]
\displaystyle\sum_{k\in N} \beta_k A_k(i) \;\le\; 0, & \forall i\in \mathcal{M}^*, \\[1em]
\beta_k \;\ge\; 0, & \forall k\in N^*.
\end{cases}
\tag{2.2}
\]

An \(n+1\)-tuple \(\beta = (\beta_0,\beta_1,\ldots,\beta_n)\) is a program of \((A,\mathcal{M})'\) if \(\beta_0=1\) and it satisfies (2.2). It is a \textbf{minimal program} if, besides, it gives the smallest value to the function $\gamma$.

Next we introduce a very important lemma, known in the literature of linear programming as complementary slackness Lemma \cite{MS}.

{ \lem\label{LP} \cite[Theorem 3.4]{D2}
\begin{enumerate}
 \item The problems $(A,\mathcal{M})$ and $(A,\mathcal{M})'$ admit at least one extremal program (i.e. a maximal and a minimal program, respectively). Each pair of programs $b$ of $(A,\mathcal{M})$ and $\beta$ of $(A,\mathcal{M})'$ satisfies $g \le \gamma$. Moreover, the extremal values of $g$ and $\gamma$ are equal.
	
\item For each pair $(b,\beta)$ of extremal programs, the following two sets of equations hold:
\begin{align}
\beta_k \left( \sum_{i \in \mathcal{M}} b_i A_k(i) \right) &= 0, \quad \forall k \in N^*, \tag{2.3} \\
b_i \left( \sum_{k \in N} \beta_k A_k(i) \right) &= 0, \quad \forall i \in \mathcal{M}^*. \tag{2.4}
\end{align}
Conversely, if a pair $(b,\beta)$ of programs satisfies $(2.3)$ and $(2.4)$, then it is a pair of extremal programs.
\end{enumerate}
}

Let \(Y \subseteq X\) be a nonempty subset and let \(R = \{R_0, R_1, \ldots, R_n\}\) be a set of relations on \(X\).
The \textbf{inner distribution} of \(Y\) with respect to \(R\) is the \((n+1)\)-tuple \(a = (a_0, a_1, \ldots, a_n)\) defined by 
\[
a_i = \frac{1}{|Y|}\,|R_i \cap Y^2|,
\] 
i.e., \(a_i\) is the average number of points of \(Y\) that are \(i\)-th associates of a fixed point of \(Y\).
Clearly, 
\[
a_0 = 1; ~\sum_{i=0}^{n} a_i = |Y|;~ a_i = a_j, ~\text{if}~ R_i^T=R_j.
\]

\section{Codes of given degree}
In this section, we adapt the proof of \cite[Theorem 5.20]{D2} to the multivariate setting. Let \(\mathfrak X=(X,\{R_\alpha\}_{\alpha\in I})\) be a multivariate \(Q\)-polynomial association scheme satisfying Hypothesis \(F_{\mathrm{deg}}\). We retain the notation \(E_\gamma\), \(\mu_\gamma\), \(z_\alpha\), and \(\Phi_\gamma\) from Section~2.2.

\subsection{Characteristic matrices}
For each $\gamma \in \mathcal{D}^*$ we construct a characteristic matrix $H_\gamma$ as follows: 
choose a real orthonormal basis of the column space of $E_\gamma$ and let $H_\gamma$ be the $|X| \times \mu_\gamma$ matrix whose columns are $\sqrt{|X|}$ times these basis vectors. Then \(H_\gamma H_\gamma^\top = |X| E_\gamma\) yields
\[
H_\gamma H_\gamma^\top = \bigl[ \Phi_\gamma(z_{\alpha(a,b)}) \bigr]_{a,b \in X}, \tag{3.1}
\]
where \(\alpha(a,b)\) denotes the unique index such that $(a,b) \in R_{\alpha(a,b)}$.
\subsection{Code \(Y\) and its parameters}
Let $Y\subseteq X$ be a code, and set $M=|Y|$. Define
\[
S = \{\alpha \in \I \setminus \{0\} : R_\alpha \cap Y^2 \neq \varnothing\},
\]
the set of non-identity relations that actually occur among pairs of elements of \(Y\), and set \(r = |S|\), the {\bf degree} of $Y$.
For every \(\beta \in S\) choose an affine functional \(L_\beta : \R^\ell \to \R\) satisfying
\[
L_\beta(z_\beta) = 1,~L_\beta(z_0) = 0.
\]
A natural choice is \(L_\beta(x) = \frac{\langle x-z_0, z_\beta-z_0 \rangle}{\|z_\beta-z_0\|^2}\).
\subsection{Annihilator polynomial}
Define the annihilator polynomial of \(Y\) by
\[
A(x) = M \prod_{\beta \in S} \bigl(1 - L_\beta(x)\bigr).
\]
Clearly \(A(z_0) = M\) and \(A(z_\beta) = 0\) for every \(\beta \in S\); its total degree is \(r\).
Let \(Z=\{z_\alpha:\alpha\in \mathcal{I}\}\) be the spectral node set, we only use the restriction of the annihilator polynomial \(A\) to \(Z\). Since the functions \(\{\Phi_\gamma|_Z:\gamma\in \mathcal{D^*}\}\) form a basis of the function space on \(Z\), there exist unique coefficients \(a_\gamma\) such that, 
\[
A(z_\alpha)= \sum_{\gamma\in \mathcal{D^*}}a_\gamma\Phi_\gamma(z_\alpha),~ \text{for every } \alpha\in \mathcal{I}.
\]
Since \(\deg A=r\), hypothesis \(F_{\mathrm{deg}}\) implies
\(A(z_\alpha)=\sum_{|\gamma|\leq r}a_\gamma\Phi_\gamma(z_\alpha), \alpha\in \mathcal{I}\).
Therefore
\[
A|_Z=\sum_{\gamma\in C}a_\gamma\Phi_\gamma|_Z,~\text{where}~ C=\{\gamma\in \mathcal{D^*}:|\gamma|\leq r\}.
\]
\subsection{Construction of the matrices \(G_r\) and \(\Delta\)}
For each \(\gamma \in C\) let \(H_\gamma(Y)\) be the submatrix of \(H_\gamma\) consisting of the rows indexed by \(Y\), it is an \(M \times \mu_\gamma\) matrix.
Concatenate these matrices horizontally:
\[
G_r = \bigl[\, H_\gamma(Y) \,\bigr]_{\gamma \in C}
~(\text{size }~M \times M(r)),
\]
where \(M(r) = \sum_{\gamma \in C} \mu_\gamma\). Define a block diagonal matrix \(\Delta\) of order \(M(r)\) by
\[
\Delta = \bigoplus_{\gamma \in C} a_\gamma I_{\mu_\gamma},
\]
i.e. each block is the scalar \(a_\gamma\) times the identity matrix of size \(\mu_\gamma\).

{\lem\label{Four}
Assume Hypothesis \(F_{\mathrm{deg}}\), with the notation above we have
\[
G_r \Delta G_r^\top = M I_M.
\]
Consequently \(M \le M(r) = \sum_{\gamma \in C} \mu_\gamma\).
}
\begin{proof}
By the expansion of the node function \(A|_Z\), for all \(a,b\in Y\),
\[
A(z_{\alpha(a,b)})=\sum_{\gamma\in C}a_\gamma\Phi_\gamma(z_{\alpha(a,b)}).
\]
Using the characteristic matrices, this gives
\[
\sum_{\gamma\in C}a_\gamma H_\gamma(Y)H_\gamma(Y)^T=
\left(A(z_{\alpha(a,b)})\right)_{a,b\in Y}.
\]
The diagonal entries on the right-hand side are \(A(z_0)=M\), while the
off-diagonal entries are zero, since \(\alpha(a,b)\in S\) whenever \(a\neq b\).
Hence \(G_r\Delta G_r^T=M I_M\). Taking ranks, we obtain
\[
M=\operatorname{rank}(M I_M)\leq \operatorname{rank}(G_r)\leq M(r).
\]
Therefore \(M\leq M(r)\).
\end{proof}

{\thm\label{Thm4.3}
Assume Hypothesis \(F_{\mathrm{deg}}\). If \(M=M(r)\), then
\begin{itemize}
\item[(1)]  The square matrix \(G_r\) of order \(M\) satisfies \(G_r G_r^\top = M I_M\).
\item[(2)]  As functions on the spectral node set \(Z\),
\[
A|_Z=Q_r|_Z,~Q_r(x)=\sum_{\gamma\in C}\Phi_\gamma(x),
\]
\(Q_r\) is called the degree Wilson polynomial, and consequently every \(z_\beta\) \((\beta \in S)\) is a zero of \(Q_r\).
\end{itemize}
}
\begin{proof} \noindent
\begin{itemize}
\item[(1)] By Lemma \ref{Four}, we have \(G_r \Delta G_r^\top = M I_M\), where \(\Delta = \bigoplus_{\gamma \in C} a_\gamma I_{\mu_\gamma}\) and \(a_\gamma\) are those in the expansion of the node function \(A|_Z\) in the \(Q\)-eigenfunction basis.
Since $M=M(r)$, $G_r$ is an $M\times M$ matrix and the right-hand side is nonsingular, hence $G_r$ is invertible. Multiplying on the left by $G_r^{-1}$ and on the right by $(G_r^{\top})^{-1}$ gives
\[
\Delta = M(G_r^{\top}G_r)^{-1}~\Longrightarrow~G_r^{\top}G_r = M\Delta^{-1}. \tag{3.2}
\]
Set $B=G_r^{\top}G_r$, then $B=M\Delta^{-1}$.
Since \(G_r\) is real and nonsingular, \(B=G_r^TG_r\) is positive definite,
hence all $a_\gamma>0$.
We compute $\operatorname{tr}(B)$ in two ways. First,
\[
\operatorname{tr}(B)=\operatorname{tr}(G_r^{\top}G_r)=\operatorname{tr}(G_rG_r^{\top}).
\]
Now $G_rG_r^{\top}=\sum_{\gamma\in C}H_\gamma(Y)H_\gamma(Y)^{\top}$, By \cite[Chapter 21, Theorem 2, p.654]{MS} and $(3.1)$, the diagonal entries of $H_\gamma(Y)H_\gamma(Y)^{\top}$ are $\Phi_\gamma(z_0)=\mu_\gamma$, hence
\[
(G_rG_r^{\top})_{x,x}=\sum_{\gamma\in C}\mu_\gamma = M(r)=M,~\forall x\in Y.
\]
Thus $\operatorname{tr}(G_rG_r^{\top})=M\cdot|Y|=M^2$, and consequently
\[
\operatorname{tr}(B)=M^2. \tag{3.3}
\]
On the other hand, using $B=M\Delta^{-1}$ and $\Delta^{-1}=\bigoplus_{\gamma\in C}a_\gamma^{-1}I_{\mu_\gamma}$, we obtain
\[
\operatorname{tr}(B)=M\sum_{\gamma\in C}\frac{\mu_\gamma}{a_\gamma}. \tag{3.4}
\]
By $(3.3)$ and $(3.4)$ we have $\sum\mu_\gamma/a_\gamma = M$.\\
Evaluating \( A|_Z=\sum_{\gamma\in C}a_\gamma\Phi_\gamma|_Z\) at \(z_0\), and using \(\Phi_\gamma(z_0)=\mu_\gamma\), gives
\[
\sum_{\gamma\in C}a_\gamma\mu_\gamma=A(z_0)=M. \tag{3.5}
\]
Moreover, by definition
\[
\sum_{\gamma\in C}\mu_\gamma = M(r)=M. \tag{3.6}
\]
Apply the Cauchy-Schwarz inequality to the vectors $(\sqrt{a_\gamma\mu_\gamma})_{\gamma\in C}$ and $(\sqrt{\mu_\gamma/a_\gamma})_{\gamma\in C}$:
\[
\biggl(\sum_{\gamma\in C}\mu_\gamma\biggr)^{\!2}
\le \biggl(\sum_{\gamma\in C}a_\gamma\mu_\gamma\biggr)
\biggl(\sum_{\gamma\in C}\frac{\mu_\gamma}{a_\gamma}\biggr).
\]
Substituting $(3.5)$, $(3.6)$ and the value $\sum\mu_\gamma/a_\gamma = M$ gives
\(M^2 \le M\cdot M = M^2\). Hence equality holds, which forces equality in the Cauchy-Schwarz inequality. Equality occurs iff there exists a constant $\lambda$ such that $\sqrt{a_\gamma\mu_\gamma}=\lambda\sqrt{\mu_\gamma/a_\gamma}$ for all $\gamma\in C$, i.e. $a_\gamma=\lambda$ (constant). Inserting this into $(3.5)$ yields $\lambda\sum\mu_\gamma=M$, and using $(3.6)$ we get $\lambda=1$. Therefore
\(a_\gamma=1,~\forall~\gamma\in C\). Thus $\Delta=I$, and substituting back into Lemma \ref{Four} gives $G_rG_r^{\top}=MI_M$,
\item[(2)] By the conclusion above, $a_\gamma=1$ for every $\gamma\in C$,
so on the spectral node set \(Z\),
\[
A|_Z=\sum_{\gamma\in C}\Phi_\gamma|_Z=Q_r|_Z.
\]
Using \(A(z_\beta) = 0\) for every \(\beta \in S\), therefore \(Q_r(z_\beta) = 0\).
\end{itemize}
\end{proof}

\section{Few distance codes}
\subsection{Distance degree of a code}
Let \(\mathfrak{X} = (X, \{R_\alpha\}_{\alpha \in \mathcal{I}})\) be an \(\ell\)-variate \(Q\)-polynomial association scheme, with notation as in Section~3.
Let \(d:\mathcal{I}\to\mathbb{R}_{\geq 0}\) be the class-distance function induced by the quasi-distance of the WMAS, so that \( \forall (x,y)\in R_\alpha,~d(x,y)=d(\alpha)\).
For a subset \(Y\subseteq X\), let \(M = |Y|\),
\[
S = \{\alpha\in\mathcal{I}\setminus\{0\}: R_\alpha\cap Y^2 \neq \varnothing\},~
D = \{d(\alpha):\alpha\in S\},~s = |D|.
\]
We call \(s\) the \textbf{distance degree} of the code \(Y\).
\subsection{Construction of the annihilator polynomial}
Assume that there exists a real affine function
\(f:\mathbb{R}^{\ell}\longrightarrow \mathbb{R}\) such that \(f(z_\alpha)=d(\alpha)\) for every \(\alpha\in\I\).
Form the univariate polynomial
\[
P(t)=\prod_{d\in D}\left(1-\frac{t}{d}\right),
\]
which satisfies \(P(0)=1\) and \(P(d)=0\) for all \(d\in D\).
Define the annihilator polynomial of \(Y\) by
\[
A(x)=M\,P(f(x))=M\prod_{d\in D}\left(1-\frac{f(x)}{d}\right).
\]
Then \(A(z_0)=M\), \(A(z_\beta)=0\) for every \(\beta\in S\), and \(\deg A = s\).
\subsection{Spectral expansion and matrix construction}
We regard the annihilator polynomial \(A\) through its restriction to the spectral node set \(Z=\{z_\alpha:\alpha\in \mathcal{I}\}\). Since \(\deg A=s\), Hypothesis \(F_{\mathrm{deg}}\) gives
\[
A|_Z=\sum_{\gamma\in C}a_\gamma\Phi_\gamma|_Z,~C=\{\gamma\in \mathcal{D^*}:|\gamma|\leq s\}.
\]
Set \(M(s)=\sum_{\gamma\in C}\mu_\gamma\).
For each \(\gamma\in C\), let \(H_\gamma(Y)\) be the submatrix of the characteristic matrix \(H_\gamma\) consisting of the rows indexed by \(Y\), it is an \(M\times\mu_\gamma\) matrix. Concatenate these matrices horizontally:
\[
G_s = [H_\gamma(Y)]_{\gamma\in C}~(\text{size }~ M\times M(s)).
\]
Define a block diagonal matrix \(\Delta\) of order \(M(s)\) by
\[
\Delta = \bigoplus_{\gamma\in C} a_\gamma I_{\mu_\gamma}.
\]

{ \lem\label{Four2}
Assume Hypothesis \(F_{\mathrm{deg}}\), with the notation above we have
\[
G_s\Delta G_s^\top = M I_M.
\]
Consequently \(M \le M(s)\).
}

\begin{proof}
The proof here is similar to Lemma \ref{Four}.
\end{proof}

Analogously to Theorem \ref {Thm4.3} we obtain the following theorem.

{\thm
Assume Hypothesis \(F_{\mathrm{deg}}\). If \(M=M(s)\), then
\begin{itemize}
\item[(1)] The square matrix \(G_s\) satisfies \(G_s G_s^\top = M I_M\).
\item[(2)] As functions on \(Z\),
\[
A|_Z=Q_s|_Z,~ Q_s(x)=\sum_{\gamma\in C}\Phi_\gamma(x).
\]
\(Q_s\) is called the distance Wilson polynomial, and consequently every point \(z_\beta\) (\(\beta\in S\)) is a zero of \(Q_s\).
\end{itemize}
}

\subsection{Additive codes}
We now specialize to additive codes in translation association schemes. Let $X$ be a finite abelian group and let $Y$ be a subgroup of $X$. Following Delsarte~\cite{D2}, define the dual of $Y$ in $X$ by
\[
Y^\perp=\left\{x^\circ\in X:\chi(x,x^\circ)=1\text{ for all }x\in Y\right\},
\]
where $\chi$ is the character pairing of $X$. We shall use the following standard duality result:
\[
(Y^\perp)^\perp=Y,~|Y|\,|Y^\perp|=|X|.
\]

Recall the classical definition of the \textbf{external distance} \(s'\) introduced by Delsarte \cite{D}: For a code \(Y\) under the Hamming metric, the external distance \(s'\)  is  the number of nonzero nontrivial coefficients in the MacWilliams transform of the distance distribution of the code. In the case where \(Y\) is linear, \(s'\) reduces to the number of non-zero weights in the orthogonal dual code \(Y^\perp\). The external distance serves as an upper bound on the covering radius \cite{S1} of \(Y\),

In WMAS, we define the \textbf{dispersion function} \cite{S3}, denoted by \(\Pi(e)\), as the number of elements \(i\) satisfying \(0 \leq d(i) \leq e\), and \(e\) is the {\bf packing radius} \cite{S1}. In metric schemes, this simplifies to \(\Pi(e) = e + 1\).

We then introduce a \textbf{generalized external distance}, denoted \(\mu\), which is designed to extend the concept to less regular metrics (such as the Lee metric or modular distance) \cite{S1}. This new parameter is defined as:
\[
\mu = s' - \Pi(e) + e + 1,
\]
where \(s'\) is the classical external distance and $e$ is the packing radius. In metric schemes, since \(\Pi(e) = e + 1\), the formula reduces to \(\mu = s'\), recovering the standard definition of external distance in such settings.
Recall \cite{S2} that a distance is {\bf graphical} if for  every pair $x,y$ with $d(x,y)>1,$  there is a $z$ such $d(x,z)=d(x,y)-1,$ and $d(z,y)=1.$

{\rem A translation association scheme is called {\bf self-dual} if it is isomorphic to its dual in the sense of \cite{D2}. See the appendix for properties and examples of self-dual translation association schemes. Throughout this subsection, we fix a self-dual identification inducing a bijection \(\sigma:\mathcal{I}\to\mathcal D^*\), with \(\sigma(0)=0\), then \(v_\alpha=\mu_{\sigma(\alpha)},~\alpha\in\mathcal{I}\).
We use the transported distance on the dual index set, defined by \(d_{\mathcal D^*}(\gamma) :=d_\mathcal{I}\bigl(\sigma^{-1}(\gamma)\bigr),~\gamma\in\mathcal D^*\).
}

{\thm\label{external}
Let \((X,\mathcal R)\) be a self-dual translation WMAS with relation index distance \(d_\mathcal{I}\), and fix the self-dual identification as in Remark~1. Assume that the induced quasi-distance on \(X\) is graphical. Let \(Y\le X\) have relation degree \(r=|S(Y)|\), as defined in Section~3, and let \(e'\) be the packing radius of \(Y^\perp\), then
\[
|Y|\le\sum_{\substack{\alpha\in \mathcal{I}\\d_\mathcal{I}(\alpha)\le\mu'}}v_\alpha,~\mu'=r-\Pi(e')+e'+1,
\]
where \(\mu'\) is the generalized external distance of \(Y^\perp\).
}
\begin{proof}
Let \(Y^\perp\) be the dual of \(Y\). By the MacWilliams duality for additive codes in translation schemes, \(s'(Y^\perp)=r\), hence its generalized external distance is \(\mu'=r-\Pi(e')+e'+1\). Then its covering radius \(\rho(Y^\perp)\) is at most \( \mu'\) by \cite[Theorem 3]{S1}. The sphere covering bound yields then
\[
|X| \leq |Y^\perp| \sum_{\substack{\alpha\in \mathcal{I}\\d_\mathcal{I}(\alpha)\le\mu'}}v_\alpha.
\]
The result follows from the identity \(|X| = |Y| \cdot |Y^\perp|\), which implies
\[
|Y| = \frac{|X|}{|Y^\perp|} \leq \sum_{\substack{\alpha\in \mathcal{I}\\d_\mathcal{I}(\alpha)\le\mu'}}v_\alpha.
\]
\end{proof}

Call a code $Y$ with $s$ distances satisfying \(|Y| = M(s)\) {\bf Generalized Hadamard code} (GH for short). We next consider additive codes $Y$ attaining the relation-degree bound \(|Y|=M(r)\). We show that their duals are perfect in the self-dual translation scheme setting.

{\coro
With the assumptions above, suppose further that the distance function dominates the total degree of the primitive idempotent indices i.e., \(d_{\mathcal D^*}(\gamma)\geq |\gamma|, ~\gamma\in\mathcal D^*\), then \(|Y| \leq M(\mu ')\). If $Y\subsetneq X$ and $|Y|=M(r)$, then \(\Pi(e')=e'+1\) and $Y^\perp$ is a perfect code with covering radius $e'$.
}
\begin{proof}
Let \(\sigma:\mathcal{I}\to\mathcal D^*\) be the self-dual identification fixed in Remark~1, then \(v_\alpha=\mu_{\sigma(\alpha)}\) and \(d_\mathcal{I}(\alpha)=d_{\mathcal D^*}(\sigma(\alpha))\) for every \(\alpha\in \mathcal{I}\). By the assumption \(d_{\mathcal D^*}(\gamma)\ge|\gamma|\), we have
\[
\{\gamma\in\mathcal D^*:d_{\mathcal D^*}(\gamma)\le\mu'\}
\subseteq
\{\gamma\in\mathcal D^*:|\gamma|\le\mu'\}.
\]
Consequently, Theorem~3 and \(|Y|\,|Y^\perp|=|X|\) give
\[
\frac{|X|}{|Y^\perp|}=|Y|
\leq \sum_{\substack{\alpha\in\mathcal{I}\\d_\mathcal{I}(\alpha)\leq\mu'}}v_\alpha
=\sum_{\substack{\gamma\in\mathcal D^*\\d_{\mathcal{D^*}}(\gamma)\le\mu'}}\mu_\gamma\\
\leq \sum_{\substack{\gamma\in\mathcal D^*\\|\gamma|\leq\mu'}}\mu_\gamma
=M(\mu')=M\bigl(r-\Pi(e')+e'+1\bigr)\leq M(r).
\]
If $|Y|=M(r)$ then $M(\mu')=M(r)$. Since \(M(r)=|Y|<|X|\), $r$ lies in the strictly increasing range of $M(\cdot)$. Hence $\mu'=r$, yielding \(\Pi(e')=e'+1\).

Moreover, equality holds throughout the preceding chain, so
\[
|Y^\perp|\sum_{\substack{\alpha\in\mathcal{I}\\ d_\mathcal{I}(\alpha)\leq r}}v_\alpha=|X|.
\]
Since Theorem \ref{external} gives \(\rho(Y^\perp)\leq r\), the radius-\(r\) balls centered at \(Y^\perp\) partition \(X\). Hence \(Y^\perp\) is \(r\)-perfect. Therefore,
\(r\leq e'\leq \rho(Y^\perp)\leq r\), and thus \(e'=\rho(Y^\perp)=r\).
\end{proof}

{\rem For a WMAS, define its \textbf{range of metricity} as the largest integer $e$ less than the class number such that $\Pi(e)=e+1$. We have computed it for several WMAS, which are given for the usual non-degenerate parameter settings of the corresponding schemes. If the scheme is not metric this parameter is usually very small.
\begin{itemize}
\item In the \underline{Lee scheme $L(n,q)$}, if $n\geq 2$ and $q\geq 4$, then $\Pi(2)=4>3$, and therefore the range is $1$. Hence, under the assumptions of Corollary~1, if an additive subgroup attains the relation-degree bound \(|Y|=M(r)\), then its perfect dual can have covering radius at most $1$.

\item In the \underline{homogeneous metric scheme} over \(\mathbb{Z}_{2^k}\) $($\(k\ge 3\)$)$, relations are indexed by triples \((\pi_U,\pi_V,\pi_S)\) with quasi-distance \(d=\pi_U+\pi_V+2\pi_S\). One finds \(\Pi(0)=1\) and \(\Pi(1)=3>2\), therefore the range is \(0\) . 

\item In the \underline{\(q\)-ary Johnson scheme} \(J_q(w,n)\) $($\(q\ge 3\), \(0<w\le n/2\)$)$, relations are indexed by ordered pairs \((i,j)\) with \(i\ge j\) and distance \(d=i+j\). For \(w\geq2\), \(\Pi(0)=1, \Pi(1)=2, \Pi(2)=4\), while for \(w=1\) we have \(\Pi(2)=3\). In either case, the range is \(1\).

\item In the \underline{sum rank scheme} with \(t\) blocks, when \(t=1\) $($rank metric$)$ we have \(\Pi(e)=e+1\) for all \(e\), the range is the class number. For \(t\ge 2\), as long as \(e\le\min_j n_j\), \(\Pi(e)=\binom{e+t}{t}\). Then \(\Pi(0)=1\) and \(\Pi(1)=t+1>2\), hence the range is \(0\).

\item In the \underline{mixed alphabet scheme}, for \(k=1\) $($classical Hamming scheme$)$ the range is the class number. For \(k\ge 2\) and \(e\le\min_j n_j\) we have \(\Pi(e)=\binom{e+k}{k}\), giving \(\Pi(0)=1\) and \(\Pi(1)=k+1>2\), thus the range is \(0\).

\item In the \underline{NRT scheme}, for \(r=1\) $($Hamming scheme$)$, the range is the class number. For $r\geq 2$ and $n\geq 2$, a direct count shows $\Pi(0)=1$, $\Pi(1)=2$, $\Pi(2)=4>3$, therefore the range is $1$.
\end{itemize}
}

\subsection{Perfect codes}
We next examine several codes that are duals  of perfect codes and compare their
sizes with the quantity \(M(s)=\sum_{|\gamma|\leq s}\mu_\gamma\), where \(s\) denotes the number of distinct non-zero distances. In particular, we determine when equality \(|C|=M(s)\) occurs in these examples. When it is not the case, it shows that a reciprocal of Corollary 1 does not hold when relation is replaced by distance.
\subsubsection{1-perfect Lee codes}
Let \(q = 2n+1\) and let \(C_n\) be the code over \(\mathbb{Z}_q\) generated by \([1,2,\dots,n]\).
Its size is \(M = q\).
In the length-$n$ Lee scheme each coordinate takes Lee values \(0,1,\dots,n\).
Primitive idempotents are indexed by multi-indices \(\gamma = (\gamma_1,\dots,\gamma_n)\) with \(\gamma_i\in\{0,1,\dots,n\}\); their multiplicities satisfy \(\mu_\gamma = \prod_{i=1}^n \mu(\gamma_i)\), where \(\mu(0)=1\) and \(\mu(k)=2\) for \(1\le k\le n\).

Magma experiments suggest that the number of distinct non-zero Lee weights of \(C_n\) is \(s = \tau(q)-1\), where \(\tau(q)\) denotes the number of positive divisors of \(q\).
We can prove this conjecture in the following two cases.
{ \prop \label{qprime}
Let $q = 2n+1$ be an odd prime and let $C$ be the linear code over $\mathbb{Z}_q$
generated by $[1,2,\dots,n]$. Then all non-zero codewords of $C$ have the same
Lee weight, hence the number of distinct non-zero Lee weights is $1$.
}
\begin{proof}
Write $\mathbb{Z}_q = \{0\} \cup [n] \cup (-[n])$ with $[n] = \{1,2,\dots,n\}$.
For $x \in \mathbb{Z}_q$ define
$|x| = \min(x \bmod q,\; q - (x \bmod q)) \in \{0,1,\dots,n\}$.
Every non-zero codeword has the form
$c_t = (t,2t,\dots,nt) \bmod q$ for some $t \in \mathbb{Z}_q^*$, and its Lee weight is
\(W_L(c_t) = \sum_{a=1}^{n} |at|\).
	
Fix $t$ and consider the map $f_t : [n] \to [n]$ defined by $f_t(a) = |at|$.
We show that $f_t$ is injective.
If $f_t(a)=f_t(b)$, then $|at|=|bt|$, which means
\[
at \equiv bt \pmod q ~ \text{or}~ at \equiv -bt \pmod q.
\]
Because $q$ is prime and $t \not\equiv 0 \pmod q$, we obtain
\[
a \equiv b \pmod q ~ \text{or}~ a \equiv -b \pmod q.
\]
The second case would give $q \mid (a+b)$. But $a,b \in [n]$ imply $2 \le a+b \le 2n < 2n+1 = q$, a contradiction. Hence $a \equiv b \pmod q$, and since $a,b \in [n]$, we conclude $a=b$. Thus $f_t$ is injective.
As $[n]$ is finite, $f_t$ is also surjective, i.e.\ a permutation of $[n]$.
Therefore, for every $t \in \mathbb{Z}_q^*$,
\[
W_L(c_t) = \sum_{a=1}^{n} f_t(a) = \sum_{x=1}^{n} x
= 1+2+\dots+n = \frac{n(n+1)}{2},
\]
which is independent of $t$.  Consequently, all non-zero codewords share the
same Lee weight, so there is exactly one distinct non-zero Lee weight.
\end{proof}

We now give a proof for the general case.
{ \prop
Let $q=2n+1$ be an odd integer, and let $C$ be the linear code over $\mathbb{Z}_q$
generated by $[1,2,\dots,n]$, then the number of distinct non-zero Lee weights is $\tau(q)-1$, where $\tau(q)$ denotes the number of positive divisors of $q$.
}
\begin{proof}
Every non-zero codeword has the form $c_t = (t,2t,\dots,nt) \bmod q$ for some $t\in\{1,\dots,q-1\}$.  Set $g = \gcd(t,q)$ and write $t = g t'$, $q = g m$, then $\gcd(t',m)=1$ and both $g,m$ are odd. For any $k\in\{1,\dots,n\}$,
$kt \bmod q = g\,(kt'\bmod m)$, and the Lee weight of the $k$-th coordinate is
\[
\|kt\|_q = \min(kt\bmod q,\; q-kt\bmod q)
= g \min(kt'\bmod m,\; m-kt'\bmod m)
= g\,\|kt'\|_m .
\]
Summing over $k$ we obtain
\begin{equation}
W(c_t) = g \sum_{k=1}^{n} \|kt'\|_m, n = \frac{q-1}{2}= \frac{gm-1}{2} \tag{4.1}.
\end{equation}
Write $g = 2h+1$ with $h = (g-1)/2$, so that $n = hm + (m-1)/2$.  Partition the summation range into $h$ complete blocks of length $m$ and a final incomplete block of length $(m-1)/2$.

In a complete block the index $k$ runs through a full residue system modulo $m$.
Because $\gcd(t',m)=1$, multiplication by $t'$ permutes $\mathbb{Z}_m$, hence
$kt'\bmod m$ takes each value $0,1,\dots,m-1$ exactly once.  The Lee weight $\|u\|_m$ is zero for $u=0$ and satisfies $\|u\|_m = \|m-u\|_m$ for $u\neq0$. Since $m$ is odd, the set $\{\|u\|_m : u=0,\dots,m-1\}$ consists of $0$ together with two copies of each integer $1,2,\dots,(m-1)/2$. Thus a complete block contributes
\[
\sum_{u=0}^{m-1}\|u\|_m = 2\sum_{i=1}^{(m-1)/2} i = \frac{m^2-1}{4}.
\]

For the incomplete block we have $k = hm + j$ with $j = 1,\dots,(m-1)/2$, and
$kt'\equiv jt'\pmod m$.  The map $j \mapsto \|jt'\bmod m\|_m$ is a permutation of
$\{1,2,\dots,(m-1)/2\}$ by the same argument as in Proposition \ref {qprime}. Therefore the incomplete block contributes
\[
\sum_{j=1}^{(m-1)/2} \|jt'\bmod m\|_m = \sum_{i=1}^{(m-1)/2} i = \frac{m^2-1}{8}.
\]

Adding the $h$ complete blocks and the incomplete block we get
\[
\sum_{k=1}^{n} \|kt'\bmod m\|_m
= \frac{g-1}{2}\cdot\frac{m^2-1}{4} + \frac{m^2-1}{8}
= \frac{g(m^2-1)}{8}.
\]
Inserting this into $(4.1)$ yields
\[
W(c_t) = g\cdot\frac{g(m^2-1)}{8}
= \frac{g^2(m^2-1)}{8}
= \frac{q^2-g^2}{8}.
\]

The formula shows that the Lee weight of $c_t$ depends only on $g=\gcd(t,q)$. The function $g \mapsto (q^2-g^2)/8$ is strictly decreasing for positive $g$, so distinct divisors $g$ give distinct weights. Hence the non-zero Lee weights are precisely the values $\frac{q^2-g^2}{8}$ with $g$ ranging over all divisors of $q$ except $q$ itself.  Their number is $\tau(q)-1$.
\end{proof}

Together with the prime case, we obtain uniformly that for every odd $q = 2n+1$,
the number of distinct non-zero Lee weights of $C_n$ equals $s = \tau(q)-1$,
where $\tau(q)$ is the number of positive divisors of $q$. And we have the following situation.
\begin{itemize}
\item If \(q\) is prime, then \(\tau(q)=2\) and \(s=1\).
The multi-indices with \(|\gamma|\le 1\) are the zero vector and the \(n\) unit vectors.
Their multiplicities are \(1\) and \(2\), respectively, whence
\(M(1) = 1 + 2n = q = M\), i.e., equality holds.
\item If \(q\) is composite, then \(s>1\).
The set of multi-indices with \(|\gamma|\le s\) strictly contains the \(1+2n\) terms above, so \(M(s) > M\).
\end{itemize}

Thus in these examples, \(M\leq M(s)\), with equality exactly when \(q\) is prime $($ i.e. when \(2n+1\) is prime$)$.

\subsubsection{Length-$2$ perfect codes in the Lee metric}
Let \(q = 2t^2+2t+1\) and let \(C_t\) be the code over \(\mathbb{Z}_q\) generated by \([1,2t+1]\). Its size is \(M = q\). Here \(d = \frac{q-1}{2} = t^2+t\) and Magma experiments suggest the following result.
{\prop
Let $t$ be a positive integer, $q = 2t^{2}+2t+1$, and let $C$ be the linear code over $\mathbb{Z}_q$ generated by $[1, 2t+1]$, then $C$ is self-dual, and the number of distinct non-zero Lee weights of $C$ is \(s = \frac{t^{2}+t}{2}\).
}
\begin{proof}
Put \(u=2t+1\), then \(u^2=(2t+1)^2=4t^2+4t+1=2q-1\equiv -1 \pmod q\). we have \((1,u)\cdot (1,u)=1+u^2\equiv 0 \pmod q\), and then \(C\subseteq C^\perp\).
Conversely, let \((x,y)\in C^\perp\), then \(x+uy\equiv 0 \pmod q\), so
\(x\equiv -uy \pmod q\). Taking \(a=x\), we obtain
\(ua\equiv ux\equiv -u^2y\equiv y \pmod q\). Thus \((x,y)=(a,ua)\in C\), we get \(C^\perp\subseteq C\), and hence \(C=C^\perp\). Therefore \(C\) is self-dual.

By \cite[Theorem 1]{GW}, $C$ is a length-\(2\) perfect Lee code: $|C| = q$, the minimum Lee distance is $2t+1$, and the Lee spheres of radius $t$ centered at the codewords tile the torus $\mathbb{Z}_q^{2}$ without overlap. A Lee sphere of radius $t$ contains exactly $q$ points, and the tiling places each point of $\mathbb{Z}_q^{2}$ in exactly one sphere. In particular, the sphere centered at the origin $O=(0,0)$ consists of the $q$ points whose $L^{1}$ (Manhattan) distance to $O$ is at most $t$, its center $O$ is labelled by the codeword $(0,0)\in C$.
Because \((2t+1)^{2} \equiv -1 \pmod q\), multiplication by $2t+1$ acts as a $90^{\circ}$ counter-clockwise rotation of the plane, hence $C$ possesses a symmetry of order $4$. The Lee weight is just the distance to the origin in the \(L^1\) distance.

We now show that distinct nonzero orbits have distinct Lee weights. For each
non-zero orbit, choose the unique representative \((x,y)\) satisfying
\(1\le x,y\le \frac{q-1}{2}\).
Since \((x,y)\in C\), we have \(y\equiv ux \pmod q\).
Thus there exists an integer \(m\) such that
\(y=ux-mq\).
The Lee weight of this representative is \(w=x+y=(u+1)x-mq\).
Suppose that two such representatives \((x,y)\) and \((x',y')\) have the same
Lee weight. Let \(y=ux-mq, y'=ux'-m'q\), then \((u+1)(x-x')=(m-m')q\).
Since \(\gcd(u+1,q)=1\), we have \(q\mid (x-x')\).
But \(1\le x,x'\le \frac{q-1}{2}\), so \(x=x'\), and then \(y=y'\). Hence the two orbits are the same, i.e. distinct nonzero rotation orbits give distinct Lee weights.

Consequently, the number of distinct non-zero Lee weights is the number of
non-zero rotation orbits: 
\[
s=\frac{q-1}{4}=\frac{2t^2+2t}{4}=\frac{t^2+t}{2}.
\]
\end{proof}
For length $2$ a multi-index is \(\gamma = (a,b)\) with \(0\le a,b\le d\) and multiplicity \(\mu(a,b)=\mu(a)\mu(b)\).  The quantity \(M(s)\) is the sum of \(\mu(a,b)\) over all pairs with \(a+b\le s\).

\begin{itemize}
\item For \(t=1\) we have \(d=2\), \(s=1\).
The admissible pairs are \((0,0)\), \((1,0)\), \((0,1)\) with multiplicities \(1,2,2\), so \(M(1) = 1+2+2 = 5 = M\).
\item For \(t=2\) we have \(d=6\), \(s=3\).
Enumerating the ten pairs with \(a+b\le 3\) yields \(M(3)=25\), while \(M=13\), hence \(M(s)>M\).
\item For larger \(t\) we have \(s \sim t^2/2\).
The number of integer pairs with \(a,b\ge0\) and \(a+b\le s\) is \(\frac{(s+1)(s+2)}{2}\sim s^2/2\).
Because \(s<d\) for all \(t\ge1\), every pair with \(a>0,b>0\) contributes multiplicity \(4\); boundary pairs (\(a=0\) or \(b=0\)) contribute \(2\) or \(1\).  Hence \(M(s) \sim 4\cdot\frac{s^2}{2}=2s^2\sim \frac{t^4}{2}\), whereas \(M \sim 2t^2\). Consequently, \(M(s)\gg M\) for sufficiently large \(t\), so equality does not occur.
\end{itemize}

These computations illustrate that equality occurs for \(t=1\), whereas
\(M(s)>M\) as \(t\) increases.

To illustrate the geometric structure of the perfect Lee code, Figure 1 shows the Lee ball ($L^1$ ball) of radius $t=2$ centered at the origin in $\mathbb{Z}_{13}^2$. It consists of all integer lattice points whose Manhattan distance from the origin is at most $t$. This ball contains exactly $q=13$ points (black dots) and lies completely in the free domain, i.e., it does not overlap with any other ball in the perfect tiling of the torus $\mathbb{Z}_{13}^2$.
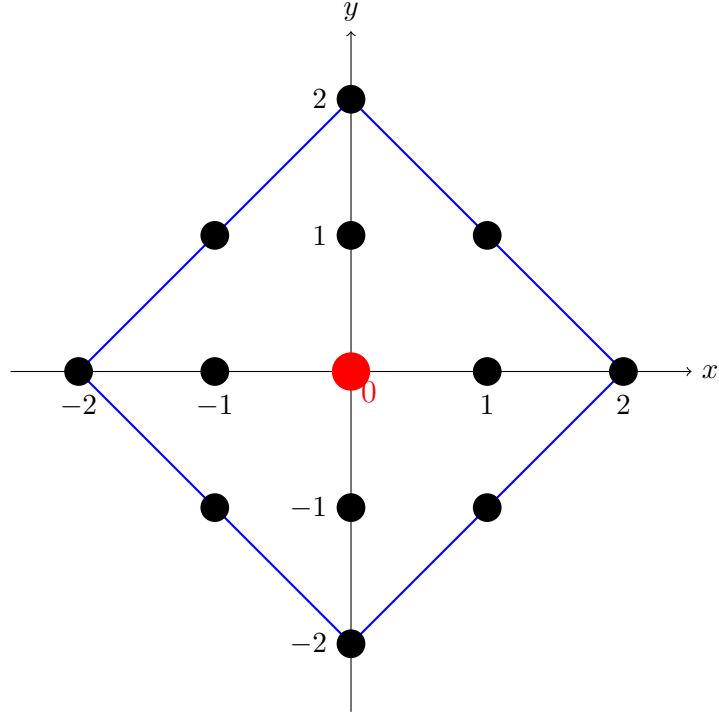
\begin{figure}[htbp]
\centering
\begin{tikzpicture}[scale=1.8]
\draw[->] (-2.5,0) -- (2.5,0) node[right] {$x$};
\draw[->] (0,-2.5) -- (0,2.5) node[above] {$y$};
\foreach \i in {-2,-1,0,1,2} {
\ifnum\i=0\relax\else
\draw (\i,0.1) -- (\i,-0.1) node[below] {$\i$};
\draw (0.1,\i) -- (-0.1,\i) node[left] {$\i$};
\fi
}
\draw[thick, blue] (2,0) -- (0,2) -- (-2,0) -- (0,-2) -- cycle;
\fill[black] (-2,0) circle (3pt);
\fill[black] (-1,-1) circle (3pt);
\fill[black] (-1,0) circle (3pt);
\fill[black] (-1,1) circle (3pt);
\fill[black] (0,-2) circle (3pt);
\fill[black] (0,-1) circle (3pt);
\fill[black] (0,0) circle (3pt);
\fill[black] (0,1) circle (3pt);
\fill[black] (0,2) circle (3pt);
\fill[black] (1,-1) circle (3pt);
\fill[black] (1,0) circle (3pt);
\fill[black] (1,1) circle (3pt);
\fill[black] (2,0) circle (3pt);
\fill[red] (0,0) circle (4pt);
\node[red, below right, font=\large] at (0,0) {$0$};
\end{tikzpicture}
\caption{Lee ball of radius $2$}
\label{fig:leeball}
\end{figure}

\subsubsection{A $1$-perfect mixed code}
Let $\mathbb{F}_8 = \mathbb{F}_2^3$ and consider the subspaces
\[
\begin{aligned}
&\mathcal{T}_1 = \{000,\;110,\;011,\;101\},\quad
\mathcal{T}_2 = \{000,\;100\},\\
&\mathcal{T}_3 = \{000,\;010\},\quad
\mathcal{T}_4 = \{000,\;001\},\quad
\mathcal{T}_5 = \{000,\;111\}.
\end{aligned}
\]
Then $\{\mathcal{T}_i\setminus\{0\}\}_{i=1}^5$ is a partition of $\mathbb{F}_8^*$.
By \cite[Theorem 7.1]{TE}, the code
\[
C = \Bigl\{(c_1,\dots,c_5)\in\prod_{i=1}^5\mathcal{T}_i :
\sum_{i=1}^5 c_i = 0 \text{ in } \mathbb{F}_8\Bigr\}
\]
is a $1$-perfect mixed code.

The size of $C$ is $|C| = 8$.
The non-zero codewords consist of six words of weight $3$ (when $c_1\neq 0$) and one word of weight $4$ (when $c_1=0$ and all later coordinates are non-zero).
Hence the set of non-zero distances is $D = \{3,4\}$ and the distance degree is $s = |D| = 2$.

The ambient space $\prod_{i=1}^5 \mathcal{T}_i$ is a direct product of five complete-graph association schemes.
Its relations are indexed by vectors $\gamma = (\gamma_1,\dots,\gamma_5)\in\{0,1\}^5$, where $\gamma_i=1$ indicates that the $i$‑th coordinate is non-zero.
The multiplicity (rank of the primitive idempotent) is \(\mu(\gamma) = 3^{\gamma_1}\),
since $\mu_1^{(1)}=3$, $\mu_1^{(i)}=1$ for $i\ge2$, and $\mu_0^{(i)}=1$ for all $i$.

For $s=2$, set
\[
M(s) = \sum_{|\gamma|\le 2}\mu(\gamma),~ |\gamma| = \sum_{i=1}^5\gamma_i.
\]
Summation by $|\gamma|$ yields:
\begin{itemize}
\item $|\gamma| = 0$: $1$ term, $\mu = 1$;
\item $|\gamma| = 1$: $\gamma_1 = 1$ ($1$ term, $\mu = 3$) and $\gamma_1 = 0$ ($\binom{4}{1} = 4$ terms, $\mu = 1$) $\rightarrow$ $3 + 4 = 7$;
\item $|\gamma| = 2$: $\gamma_1 = 1$ ($\binom{4}{1} = 4$ terms, $\mu = 3$) $\rightarrow$ $12$; $\gamma_1 = 0$ ($\binom{4}{2} = 6$ terms, $\mu = 1$) $\rightarrow$ $6$; total $18$.
\end{itemize}
Thus $M(s) = 1+7+18 = 26$, and hence \(|C|=8<26=M(s)\).

\section{ Some multivariate Q-polynomial schemes}
To illustrate the wide applicability of the bound obtained in Section $4$, we now examine several association schemes that are known to be multivariate $Q$-polynomial. For each scheme, it is enough to verify Hypothesis $F_{\deg}$ and the existence of an affine
function $f$ such that \(f(z_\alpha)=d_\mathcal{I}(\alpha),~\alpha\in\mathcal{I}\). The following observation allows us to use convenient affine parametrizations of the spectral nodes.

{\rem\label{z_a}
Let \(Z=\{z_\alpha:\alpha\in\mathcal{I}\}\) be the spectral node set. Suppose that
\(\widetilde z_\alpha=Tz_\alpha+b\), where $T\in GL_\ell(\mathbb R)$, and define
\(\widetilde\Phi_\gamma(x)=\Phi_\gamma\bigl(T^{-1}(x-b)\bigr)\), then
\[
Q_\gamma(\alpha)=\widetilde\Phi_\gamma(\widetilde z_\alpha),~\deg\widetilde\Phi_\gamma=\deg\Phi_\gamma.
\]
Hence Hypothesis $F_{\deg}$ is invariant under invertible affine reparametrizations of the spectral nodes. The existence of an affine function satisfying \(f(z_\alpha)=d_\mathcal{I}(\alpha)\) is invariant as well. Whenever such a reparametrization is used below, we suppress the tildes and continue to write $z_\alpha$ and $\Phi_\gamma$.
}

The next lemma gives a convenient criterion for Hypothesis $F_{\deg}$

{\lem\label{Fdeg}
Let \(\mathfrak{X}=(X,\{R_\alpha\}_{\alpha\in\mathcal{I}})\) be an $\ell$-variate $Q$-polynomial association scheme with dual index domain \(\mathcal D^*\subseteq\mathbb N^\ell\). Suppose that, after an invertible affine reparametrization of the spectral nodes, there exists a bijection \(\tau:\mathcal{I}\longrightarrow\mathcal D^*\) such that
\(z_\alpha=\tau(\alpha),~\alpha\in\mathcal{I}\). If \(\deg\Phi_\gamma\le |\gamma|,~\gamma\in\mathcal D^*\), then Hypothesis $F_{\deg}$ holds.
}
\begin{proof}
For $s\ge0$, set
\[
V_s=\{F|_Z:\deg F\le s\},~W_s=\operatorname{span}\{\Phi_\gamma|_Z:\gamma\in\mathcal D^*,\ |\gamma|\le s\}.
\]
Since $\deg\Phi_\gamma\le|\gamma|$, we have \(W_s\subseteq V_s\). The second eigenmatrix $Q$ is nonsingular, and \(Q_\gamma(\alpha)=\Phi_\gamma(z_\alpha)\), so the functions \(\{\Phi_\gamma|_Z:\gamma\in\mathcal D^*\}\) are linearly independent. Hence
\[
\dim W_s=\#\{\gamma\in\mathcal D^*:|\gamma|\le s\}.
\]
By Definition~6, $\mathcal D^*$ is downward closed under the componentwise order. For $\beta\in\mathbb N^\ell$, define \(N_\beta(x)=\prod_{i=1}^{\ell}\binom{x_i}{\beta_i}\).
The polynomials $N_\beta$ with $|\beta|\le s$ form a basis of the polynomials of total degree at most $s$. Moreover, if $\beta\notin\mathcal D^*$, then \(N_\beta|_Z=0\), since $N_\beta(\lambda)\ne0$ for some $\lambda\in\mathcal D^*$ would imply $\beta\le\lambda$, and hence $\beta\in\mathcal D^*$ by downward closure.
Since $N_\beta(\lambda)\neq 0$ only if $\beta\leq\lambda$, the matrix \(\bigl(N_\beta(\lambda)\bigr)_{\lambda,\beta\in\mathcal D^*}\)
can be arranged in triangular form, with diagonal entries $N_\lambda(\lambda)=1$. Hence it is nonsingular, and the functions $N_\beta|_Z$, $\beta\in\mathcal D^*$ are linearly independent, therefore
\[
\dim V_s=\#\{\beta\in\mathcal D^*:|\beta|\le s\}.
\]
Thus $\dim V_s=\dim W_s$. Together with $W_s\subseteq V_s$, this gives
\(V_s=W_s\), and then Hypothesis $F_{\deg}$ holds.
\end{proof}

For each scheme below, we first show that after an invertible affine reparametrization, the spectral nodes can be identified with the dual index domain $\mathcal D^*$, and verify that \(\deg\Phi_\gamma\leq |\gamma|,~\gamma\in\mathcal D^*\). Lemma~\ref{Fdeg} then yields Hypothesis $F_{\deg}$. We finally exhibit an affine function $f$ satisfying
\(f(z_\alpha)=d_\mathcal{I}(\alpha),~\alpha\in\mathcal{I}\).

\subsection{$q$-ary Johnson scheme}
Let $q\geq 3$ and $0<w<n$. Following \cite[Section~3.2]{Johnson}, the relations of the
$q$-ary Johnson scheme $J_q(n,w)$ are indexed by
\[
D=\{(a,b)\in\mathbb N^2:a+b\leq w,\;0\leq b\leq \min\{w,n-w\}\},
\]
and the primitive idempotents are indexed by the same set $D$. Moreover, $J_q(n,w)$ is bivariate $Q$-polynomial with respect to the graded lexicographic order \cite[Theorem~1.1]{Johnson}. The domain $D$ is componentwise downward closed.

By the second-eigenvalue formulas in \cite[(3.10)]{Johnson},
\[
\begin{aligned}
Q_{1,0}(a,b)&=\frac{n}{w}\bigl((q-2)(w-b)-(q-1)a\bigr),\\
Q_{0,1}(a,b)&=\frac{n-1}{w(n-w)}\bigl(w(n-w)-nb\bigr).
\end{aligned}
\]
Thus the degree-one spectral coordinates are an invertible affine function of $(a,b)$, since their linear part is triangular with nonzero diagonal entries. By Remark~\ref{z_a}, we may therefore take \(z_{a,b}=(a,b)\).

Since the graded lexicographic order is $L^1$-compatible and $\Phi_\gamma$ has multidegree $\gamma$, we have \(\deg\Phi_\gamma=|\gamma|\). Hence Lemma~\ref{Fdeg} implies that Hypothesis $F_{\deg}$ holds.

Finally, for a relation indexed by $(a,b)$, the Hamming distance is
\(d_\mathcal{I}(a,b)=a+2b\). Therefore, taking \(f(x_1,x_2)=x_1+2x_2\), we obtain
\(f(z_{a,b})=d_\mathcal{I}(a,b)\).

\subsection{Lee scheme}
Let $s=\lfloor q/2\rfloor$. For $n>1$, the Lee scheme on $\mathbb Z_q^n$ is the Delsarte extension of the one-coordinate Lee scheme \cite{Lee}. Its relations are indexed by the Lee composition parameters
\[
D=\left\{\alpha=(k_1,\ldots,k_s)\in\mathbb N^s:\sum_{i=1}^s k_i\le n\right\},
\]
where $k_i$ is the number of coordinates of Lee weight $i$. The set $D$ is componentwise downward closed. By the extension theorem \cite[Theorem~4.3]{B+}, the relations and primitive idempotents are both indexed by $D$, and the Lee scheme is $s$-variate $Q$-polynomial with respect to the graded lexicographic order.

Let $Q^{(1)}$ denote the second eigenmatrix of the one-coordinate Lee scheme, and put
\(k_0=n-\sum_{r=1}^s k_r\). By the Lee-number generating formula for the Delsarte extension \cite{Lee}, the degree-one $Q$-eigenvalues satisfy
\[
Q_{\epsilon_j}(\alpha)=\sum_{r=0}^s k_rQ_j^{(1)}(r)=nQ_j^{(1)}(0)+\sum_{r=1}^s\bigl(Q_j^{(1)}(r)-Q_j^{(1)}(0)\bigr)k_r .
\]
Hence the spectral coordinates are affine-linear functions of $\alpha=(k_1,\ldots,k_s)$. Their linear part is nonsingular, since subtracting the $0$-th column of $Q^{(1)}$ from each of its remaining columns shows that its determinant is $\det Q^{(1)}\neq0$.
Therefore, by Remark~\ref{z_a}, we may take \(z_\alpha=\alpha=(k_1,\ldots,k_s)\).

Since the graded lexicographic order is $L^1$-compatible and $\Phi_\gamma$ has multidegree $\gamma$, we have \(\deg\Phi_\gamma=|\gamma|\). Therefore Lemma~\ref{Fdeg} yields Hypothesis $F_{\deg}$.

Finally, the Lee distance of the relation indexed by $\alpha=(k_1,\ldots,k_s)$ is
\(d_\mathcal{I}(\alpha)=\sum_{i=1}^s i k_i\). Therefore, taking
\(f(x_1,\ldots,x_s)=\sum_{i=1}^s i x_i\), we have \(f(z_\alpha)=d_\mathcal{I}(\alpha)\).

\subsection{Mixed alphabet scheme}
Let \(X=\prod_{j=1}^{m}\mathbb Z_{p_j}^{\,n_j}\) equipped with the Hamming distance. The corresponding association scheme is the direct product \(\mathfrak{X}=\bigotimes_{j=1}^{m}H(n_j,p_j)\). By \cite[Theorem~4.1]{B+}, the relations and primitive idempotents are both indexed by \(D=\prod_{j=1}^{m}\{0,1,\ldots,n_j\}\), and $\mathfrak{X}$ is $m$-variate $Q$-polynomial. The domain $D$ is componentwise downward closed.

For
\(\alpha=(s_1,\ldots,s_m)\in D\), where $s_j$ is the Hamming distance in the $j$-th component, the degree-one $Q$-eigenvalues are \cite{Mixed}
\[
Q_{\epsilon_j}(\alpha)=(p_j-1)n_j-p_js_j,~1\leq j\leq m,
\]
by the degree-one Krawtchouk polynomial of $H(n_j,p_j)$. Hence the spectral coordinates are related to $\alpha$ by an invertible affine transformation, whose linear part is
\(\operatorname{diag}(-p_1,\ldots,-p_m)\). By Remark~\ref{z_a}, we may therefore take
\(z_\alpha=\alpha=(s_1,\ldots,s_m)\).

The second eigenvalues of the direct product are products of the Krawtchouk polynomials of the component Hamming schemes \cite{DP}. Thus, for
\(\gamma=(\gamma_1,\ldots,\gamma_m)\in D\), we may write
\[
\Phi_\gamma(x_1,\ldots,x_m)=\prod_{j=1}^{m}K_{\gamma_j}^{(p_j)}(x_j)\,
\]
and hence \(\deg\Phi_\gamma=\sum_{j=1}^{m}\gamma_j=|\gamma|\). Therefore Lemma~\ref{Fdeg} yields Hypothesis $F_{\deg}$.

Finally, \(d_\mathcal{I}(\alpha)=\sum_{j=1}^{m}s_j\). Therefore, taking \(f(x_1,\ldots,x_m)=\sum_{j=1}^{m}x_j\), we obtain \(f(z_\alpha)=d_\mathcal{I}(\alpha)\).

\subsection{NRT scheme}
Let $q$ be a prime power and let $\mathbb{F}_q^{rn}$ be partitioned into $n$ blocks of length $r$. For a vector $x$, let \(\alpha=(\lambda_1,\ldots,\lambda_r)\) be its shape, where $\lambda_i$ is the number of blocks whose rightmost nonzero entry is in position $i$. The relations of the NRT scheme are indexed by
\[
D=\Delta_{r,n}=\left\{(\lambda_1,\ldots,\lambda_r)\in\mathbb{N}^r:\sum_{i=1}^r\lambda_i\leq n\right\},
\]
which is componentwise downward closed. By \cite[Theorem~3.1]{NR}, the ordered Hamming scheme is the $n$-fold Delsarte extension of an $r$-class association scheme. Hence,
by \cite[Theorem~4.3]{B+}, its relations and primitive idempotents are both indexed by $D$, and it is an $r$-variate $Q$-polynomial scheme with respect to the graded lexicographic order.

The degree-one multivariate Krawtchouk polynomials are given by \cite[Eq.~(21)]{NR} as
\[
K_{F_i}(x)=q^{i-1}(q-1)\bigl(n-x_r-\cdots-x_{r-i+2}\bigr)-q^i x_{r-i+1},~1\leq i\leq r.
\]
Thus the degree-one spectral coordinates are affine-linear functions of the shape parameters. After reversing the order of the variables, their linear part is triangular with diagonal entries \(-q,-q^2,\ldots,-q^r\), and is therefore nonsingular. By Remark~\ref{z_a}, we may hence take \(z_\alpha=\alpha=(\lambda_1,\ldots,\lambda_r)\).

Moreover, the eigenvalues of the ordered Hamming scheme are the multivariate Krawtchouk polynomials: \(Q_\gamma(\alpha)=K_\gamma(\alpha)\), and $K_\gamma$ has total degree $|\gamma|$ \cite[Section~4-A]{NR}. Consequently, \(\deg\Phi_\gamma=|\gamma|\), and Lemma~\ref{Fdeg} yields Hypothesis $F_{\deg}$.

Finally, the NRT distance of the relation indexed by $\alpha=(\lambda_1,\ldots,\lambda_r)$ is \(d_\mathcal{I}(\alpha)=\sum_{i=1}^r i\lambda_i\). Therefore, taking \(f(x_1,\ldots,x_r)=\sum_{i=1}^r i x_i\), we have
\(f(z_\alpha)=d_\mathcal{I}(\alpha)\).

\section{Designs}
\begin{definition}
Let $(X,R)$ be an $\ell$-variable $Q$-polynomial scheme, and assume $(X,R)$ is a WMAS with distance function $d=d_{\mathcal{D}^*}: \mathcal{D}^* \to \mathbb R_{\geq 0}$.
Let $T \subseteq \mathcal{D}^*$ be arbitrary, with $0 \notin T$.
A set $Y \subseteq X$ is a \textbf{$T$-design} if its inner distribution $a$ satisfies
\[
(aQ)_i = 0 \quad \text{for all } i \in T.
\]
For a positive integer $t$, $Y$ is a \textbf{$t$-design} if it is a $T_t$-design with
\[
T_t = \{ i \in \mathcal{D}^* \mid 0 < d(i) \leq t \}.
\]
\end{definition}

In the next two subsections we study each of these two concepts of designs in turn, and give an analogue of the Rao bound in each case.
\subsection{t-designs}
Before applying the linear-programing framework of Section~2.3, we fix the following indexing convention. Let \(\mathcal{I}=\{\alpha_0,\alpha_1,\ldots,\alpha_n\}, \mathcal{D}^*=\{\gamma_0,\gamma_1,\ldots,\gamma_n\}\) with $\alpha_0=\gamma_0=0$, and identify both enumerations with the abstract index set \(N=\{0,1,\ldots,n\}\). Define \(\widehat{P}_k(i)=P_{\alpha_k}(\gamma_i),~k,i\in N\). Since \(P_0(\gamma)=1,~\gamma\in\mathcal{D}^*\) and \(P_\alpha(0)=v_\alpha>0,~\alpha\in\mathcal{I}\), the matrix \(\widehat{P}=[\widehat{P}_k(i)]_{k,i \in N}\) satisfies the assumptions imposed on the matrix $A$ in Section~2.3.

{\lem\label{L-P}
Let \(Y\subseteq X\) be nonempty,and let \(a=(a_\alpha)_{\alpha\in\mathcal{I}}\) be its inner distribution. Define \(b_i=\frac{(aQ)_{\gamma_i}}{|Y|},~i\in N\), then
\(b_0=1\), \(b_i\geq0,~i\in N\), and \(bP=\frac{|X|}{|Y|}a\).
}
\begin{proof}
Since \(Q_0(\alpha)=1, \forall \alpha\in\mathcal{I}\) and \(\sum_{\alpha\in\mathcal{I}}a_\alpha=|Y|\), we have \((aQ)_0=\sum_{\alpha\in\mathcal{I}}a_{\alpha}Q_0(\alpha)=|Y|\), and then \(b_0=1\).
Moreover, using \(|X|E_\gamma=\sum_{\alpha\in\mathcal{I}}Q_\gamma(\alpha)A_\alpha\), we obtain
\[
(aQ)_\gamma=\frac{|X|}{|Y|}\chi_Y^{T}E_\gamma\chi_Y=\frac{|X|}{|Y|}\|E_\gamma\chi_Y\|^2\geq0,
\]
thus \(b_i\geq0,~i\in N\). Finally, by \(QP=|X|I\), we have
\[
bP=\frac{1}{|Y|}aQP=\frac{|X|}{|Y|}a.
\]
\end{proof}

\begin{definition}
For $u,v\in\mathcal{D}^*$, the {\bf Krein parameters} $q_{uv}^{w}$ are defined by
\[
E_u\circ E_v=\frac{1}{|X|}\sum_{w\in\mathcal{D}^*}q_{uv}^{w}E_w.
\]
\end{definition}
Krein parameters \cite[Section~11.4]{uvp} satisfy \(q_{uv}^{0}=\delta_{uv}\mu_u\). Moreover, there has Krein product formula \(Q_u(\alpha)Q_v(\alpha)=\sum_{w\in\mathcal{D}^*}q_{uv}^{w}Q_w(\alpha),~\alpha\in\mathcal{I}\).

{\thm\label{Distance Rao bound}
Let $\mathfrak{X}=(X,\{R_\alpha\}_{\alpha\in\mathcal{I}})$ be a symmetric multivariate
$Q$-polynomial WMAS, and let $Y\subseteq X$ be a $t$-design. Assume that the dual-index distance \(d=d_{\mathcal{D}^*}:\mathcal{D}^*\longrightarrow\mathbb{R}_{\geq0}\) is compatible with the Krein support, namely,
\[
q_{uv}^{w}\neq0 ~\Longrightarrow ~d(w)\leq d(u)+d(v), ~u,v,w\in\mathcal{D}^*.
\]
Set \(e=\left\lfloor\frac{t}{2}\right\rfloor\) and \(\mu_e'
=\sum_{\substack{\gamma\in\mathcal{D}^*\\ d(\gamma)\leq e}}\mu_\gamma\),
then \(|Y|\geq\mu_e'\).
}
\begin{proof}
Define \(\Psi_e(z)=\sum_{\substack{u\in\mathcal{D}^*\\d(u)\leq e}}\Phi_u(z)\) and \(\beta_k=\left(\frac{\Psi_e(z_{\alpha_k})}{\mu_e'}\right)^2,~k\in N\).
Since \(\Phi_u(z_0)=Q_u(0)=\mu_u\), we have
\[
\Psi_e(z_0)=\sum_{\substack{u\in\mathcal{D}^*\\d(u)\leq e}}\Phi_u(z_0)=\sum_{\substack{u\in\mathcal{D}^*\\d(u)\leq e}}\mu_u=\mu_e'.
\]
Hence $\beta_0=1$, and clearly $\beta_k\geq0$ for all $k\in N$.
For $\gamma\in\mathcal{D}^*$, the identity
\(\mu_\gamma P_\alpha(\gamma)=v_\alpha Q_\gamma(\alpha)=v_\alpha\Phi_\gamma(z_\alpha)\)
gives
\[
\begin{aligned}
\sum_{k\in N}\beta_kP_{\alpha_k}(\gamma)
&=
\frac{1}{\mu_\gamma(\mu_e')^2}
\sum_{\alpha\in\mathcal{I}}
v_\alpha\Phi_\gamma(z_\alpha)\Psi_e(z_\alpha)^2.
\end{aligned}
\]
By the Krein product formula, we have
\(\Phi_u(z_\alpha)\Phi_v(z_\alpha)=\sum_{w\in\mathcal{D}^*}q_{uv}^{w}\Phi_w(z_\alpha)\),
so that
\[
\Psi_e(z_\alpha)^2=\sum_{w\in\mathcal{D}^*}c_w\Phi_w(z_\alpha),~\text{where}~
c_w=\sum_{\substack{u,v\in\mathcal{D}^*\\d(u)\leq e,\ d(v)\leq e}}q_{uv}^{w}.
\]
If $c_w\neq0$, then $q_{uv}^{w}\neq0$ for some $u,v$ satisfying $d(u),d(v)\leq e$.
By the Krein-support compatibility assumption, \(d(w)\leq d(u)+d(v)\leq2e\).
Therefore whenever \(d(w)>2e\), \(c_w=0\).
By Lemma~\ref{eigenmatrix-identities} \(\sum_{\alpha\in \mathcal{I}}v_\alpha\Phi_\gamma(z_\alpha)\Phi_w(z_\alpha)=|X|\mu_\gamma\delta_{\gamma,w}\), we obtain
\[
\sum_{k\in N}\beta_kP_{\alpha_k}(\gamma)=\frac{|X|}{(\mu_e')^2}c_\gamma.
\]
Define the index set for the linear program by
\(\mathcal{M}=\{0\}\cup\{\,i\in N^*:d(\gamma_i)>t\,\}\).
Let $i\in\mathcal{M}^*$, then \(d(\gamma_i)>t\geq 2e\). Thus $c_{\gamma_i}=0$, and then
\[
\sum_{k\in N}\beta_k\widehat{P}_k(i)=\sum_{k\in N}\beta_kP_{\alpha_k}(\gamma_i)=\frac{|X|}{(\mu_e')^2}c_{\gamma_i}=0.
\]
We verify that $\beta=(\beta_k)_{k\in N}$ is feasible for the dual program $(\widehat{P},\mathcal{M})'$.
For $\gamma=0$, using Lemma~\ref{eigenmatrix-identities} \(q_{uv}^{0}=\delta_{uv}\mu_u\),
we obtain
\[
c_0=\sum_{\substack{u,v\in\mathcal{D}^*\\d(u)\leq e,\ d(v)\leq e}}q_{uv}^{0} =\sum_{\substack{u\in\mathcal{D}^*\\d(u)\leq e}}\mu_u=\mu_e'.
\]
Therefore, $\beta$ is feasible for $(\widehat{P},\mathcal{M})'$ with objective value
\[
\sum_{k\in N}\beta_k\widehat{P}_k(0)=\frac{|X|}{(\mu_e')^2}c_0=\frac{|X|}{\mu_e'}.
\]
Let $a=(a_\alpha)_{\alpha\in\mathcal{I}}$ be the inner distribution of $Y$, and define
\(b_i=\frac{(aQ)_{\gamma_i}}{|Y|},~i\in N\). By Lemma~\ref{L-P}, \(b_0=1\);~\(b_i\geq0,~i \in\mathcal{M^*}\). If $i\notin\mathcal{M}$, then $i\neq0$ and \(0<d(\gamma_i)\leq t\).
Since $Y$ is a $t$-design, we have \((aQ)_{\gamma_i}=0\), and hence \(b_i=0\).
For every $k\in N^*$,
\[
\sum_{i\in\mathcal{M}}b_i\widehat{P}_k(i)=\sum_{i\in N}b_iP_{\alpha_k}(\gamma_i)=\frac{|X|}{|Y|}a_{\alpha_k}.
\]
Because inner distribution \(a_{\alpha_k}\geq0\), then \(\sum_{i\in\mathcal{M}}b_i\widehat{P}_k(i)\geq0\), hence $b$ is feasible for $(\widehat{P},\mathcal{M})$ with objective value
\[
\sum_{i\in\mathcal{M}}b_i=\sum_{i\in N}b_i=\sum_{i\in N}b_iP_0(\gamma_i)=(bP)_0=\frac{|X|}{|Y|}a_0=\frac{|X|}{|Y|}.
\]
By Lemma~\ref{LP}, \(\frac{|X|}{|Y|}\leq\frac{|X|}{\mu_e'}\), and therefore
\begin{equation}\label{distace Rao bound}
|Y|\geq\mu_e'.
\end{equation}
\end{proof}

We call \eqref{distace Rao bound} the {\bf distance Rao bound}, and the sum polynomial
\(\Psi_e(x)=\sum_{\substack{u\in\mathcal{D}^*\\d(u)\leq e}}\Phi_u(x)\) {\bf distance Wilson polynomial}.

\begin{definition}
A $t$-design $Y$ in a multivariate $Q$-polynomial scheme will be said to be a \textbf{distance tight t-design} of order $e$ if it satisfies the distance Rao bound, i.e.
\[
|Y|=\sum_{\substack{\gamma\in\mathcal{D}^*\\ d(\gamma)\leq e}}\mu_\gamma,
\]
with $e = \lfloor t/2\rfloor$.
\end{definition}

We now relate distance-tight designs to perfect codes in the self-dual translation setting. Fix a compatible self-dual identification
\[
\sigma:\mathcal{I}\longrightarrow\mathcal{D}^*,~v_\alpha=\mu_{\sigma(\alpha)},~d_\mathcal{I}(\alpha)=d_{\mathcal{D}^*}(\sigma(\alpha)),~\forall \alpha\in\mathcal{I}.
\]

{\prop
Let \(\mathfrak X=(X,\{R_\alpha\}_{\alpha\in \mathcal{I}})\) be a self-dual
translation multivariate \(Q\)-polynomial scheme with relation index distance \(d_I\). Fix the self-dual identification and the transported dual index distance \(d_{\mathcal D^*}\) as in Remark~1, and assume that the  quasi-distance on \(X\) is graphical. Let $Y\subseteq X$ be an additive code, then $Y$ is an $e$-perfect code if and only if $Y^\perp$ is a distance-tight $2e$-design.
}
\begin{proof}
By the generalized MacWilliams identity for translation association
schemes \cite[Section~6.2]{D2}, applied to $C=Y^\perp$, we have
\[
(a(Y^\perp)Q)_{\sigma(\alpha)}=|Y^\perp|\,a_\alpha(Y),~\alpha\in\mathcal{I}.
\]
Hence, using the distance-preserving identification,
\begin{equation}\label{design-min-distance}
\begin{aligned}
Y^\perp \text{ is a }2e\text{-design}
&\Longleftrightarrow
(a(Y^\perp)Q)_\gamma=0 ~\forall\,\gamma\in\mathcal{D}^*:0<d_{\mathcal{D}^*}(\gamma)\leq2e
\\
&\Longleftrightarrow
a_\alpha(Y)=0 ~\forall\,\alpha\in\mathcal{I}:0<d_\mathcal{I}(\alpha)\leq2e
\\
&\Longleftrightarrow
d_{\min}(Y)>2e.
\end{aligned}
\end{equation}
Assume first that $Y$ is an $e$-perfect code. Since $d$ is graphical, the radius-$e$ balls centered at distinct codewords cannot be disjoint unless \(d_{\min}(Y)>2e\). Thus by \eqref{design-min-distance}, $Y^\perp$ is a $2e$-design. Moreover, the perfect
sphere partition gives
\[
|X|=|Y||B(x,e)|=|Y|\sum_{\substack{\alpha\in\mathcal{I}\\d_\mathcal{I}(\alpha)\leq e}}v_\alpha.
\]
Using $|Y||Y^\perp|=|X|$, we obtain
\[
|Y^\perp|=\sum_{\substack{\alpha\in\mathcal{I}\\d_\mathcal{I}(\alpha)\leq e}}v_\alpha=
\sum_{\substack{\gamma\in\mathcal{D}^*\\d_{\mathcal{D}^*}(\gamma)\leq e}}\mu_\gamma.
\]
Thus $Y^\perp$ is a distance-tight $2e$-design.\\
Conversely, suppose that $Y^\perp$ is a distance-tight $2e$-design, then
\[
d_{\min}(Y)>2e,~\text{and}~|Y^\perp|=\sum_{\substack{\gamma\in\mathcal{D}^*\\
d_{\mathcal{D}^*}(\gamma)\leq e}}\mu_\gamma=\sum_{\substack{\alpha\in\mathcal{I}\\
d_\mathcal{I}(\alpha)\leq e}}v_\alpha.
\]
Since $|Y||Y^\perp|=|X|$,
\[
|Y|\sum_{\substack{\alpha\in\mathcal{I}\\d_\mathcal{I}(\alpha)\leq e}}v_\alpha=|X|.
\]
The condition $d_{\min}(Y)>2e$ implies that the radius-$e$ balls centered at distinct codewords are pairwise disjoint. Hence $Y$ is an $e$-perfect code.
\end{proof}

\subsection{ T-designs}
{\thm \label{Degree Rao bound}
Let $\mathfrak{X}=(X,\{R_\alpha\}_{\alpha\in\mathcal{I}})$ be a symmetric multivariate
$Q$-polynomial association scheme with respect to an $L^1$-compatible monomial order, and assume that Hypothesis $F_{\mathrm{deg}}$ holds for its spectral node set. Let $0\notin T\subseteq\mathcal{D}^*$, and put \(\rho^*=\max_{\gamma\in\mathcal{D}^*}|\gamma|\).
Define
\[
e(T)=\max\left\{m\in\{0,1,\ldots,\rho^*\}:\{\gamma\in\mathcal{D}^*:0<|\gamma|\leq 2m\}\subseteq T\right\}.
\]
Set
\(\mu_{e(T)}'=\sum_{\substack{\gamma\in\mathcal{D}^*\\|\gamma|\leq e(T)}}\mu_\gamma\), then every $T$-design $Y\subseteq X$ satisfies \(|Y|\geq \mu_{e(T)}'\).
}
\begin{proof}
Set \(e=e(T)\), define \(\Psi_e(z)=\sum_{\substack{u\in\mathcal{D}^*\\|u|\leq e}}\Phi_u(z)\), \(\beta_k=\left(\frac{\Psi_e(z_{\alpha_k})}{\mu_e'}\right)^2\).
Since \(\Phi_u(z_0)=\mu_u\), we have $\Psi_e(z_0)=\mu_e'$, and hence $\beta_0=1$.
Clearly $\beta_k\geq0$ for all $k\in N$. As Theorem~\ref{Distance Rao bound}, for $\gamma\in\mathcal{D}^*$, we obtain
\[
\sum_{k\in N}\beta_kP_{\alpha_k}(\gamma)=\frac{1}{\mu_\gamma(\mu_e')^2}\sum_{\alpha\in\mathcal{I}}
v_\alpha\Phi_\gamma(z_\alpha)\Psi_e(z_\alpha)^2.
\]
For each $u\in\mathcal{D}^*$, the polynomial $\Phi_u$ has multidegree $u$. Since the monomial order is $L^1$-compatible, its ordinary total degree is $|u|$, hence
\(\deg\Psi_e\leq e,~\deg\Psi_e^2\leq2e\). By Hypothesis $F_{\mathrm{deg}}$,
\(\Psi_e^2\big|_Z\in\operatorname{span}\left\{\Phi_w\big|_Z:w\in\mathcal{D}^*,\ |w|\leq2e\right\}\).
On the other hand, by the Krein product formula,
\[
\Psi_e(z_\alpha)^2=\sum_{w\in\mathcal{D}^*}c_w\Phi_w(z_\alpha),~\text{where}~c_w=\sum_{\substack{u,v\in\mathcal{D}^*\\|u|\leq e,\ |v|\leq e}}q_{uv}^{w}
\]
Since the functions $\{\Phi_w|_Z:w\in\mathcal{D}^*\}$ are linearly independent, it follows that whenever \(|w|>2e\), \(c_w=0\). Exactly as in Theorem~\ref{Distance Rao bound}, we then obtain
\[
\sum_{k\in N}\beta_kP_{\alpha_k}(\gamma)=\frac{|X|}{(\mu_e')^2}c_\gamma.
\]
Define \(\mathcal{M}=\{0\}\cup\{\,i\in N^*:\gamma_i\notin T\,\}\).
If $i\in\mathcal{M}^*$, then $\gamma_i\notin T$, and hence $|\gamma_i|>2e$. Therefore $c_{\gamma_i}=0$, so \(\sum_{k\in N}\beta_k\widehat{P}_k(i)=0,~i\in\mathcal{M}^*\).
Consequently, $\beta$ is feasible for the dual program $(\widehat{P},\mathcal{M})'$ with objective value \(\frac{|X|}{\mu_e'}\).
Let $a=(a_\alpha)_{\alpha\in\mathcal{I}}$ be the inner distribution of $Y$ and define
\(b_i=\frac{(aQ)_{\gamma_i}}{|Y|},~i\in N\). If $i\notin\mathcal{M}$, then $\gamma_i\in T$, and since $Y$ is a $T$-design, \((aQ)_{\gamma_i}=0\). Hence $b_i=0$ for $i\notin\mathcal{M}$. As Theorem~\ref{Distance Rao bound}, we have
\[
b_0=1;~b_i\geq0; ~ \sum_{i\in\mathcal{M}}b_i\widehat{P}_k(i)\geq0.
\]
Therefore, $b$ is feasible for the primal program $(\widehat{P},\mathcal{M})$ with objective value \(\frac{|X|}{|Y|}\).
By Lemma~\ref{LP}, \(\frac{|X|}{|Y|}\leq\frac{|X|}{\mu_e'}\), and therefore
\begin{equation}\label{degree Rao bound}
|Y|\geq\mu_e'=\mu_{e(T)}'.
\end{equation}
\end{proof}

We call \eqref{degree Rao bound} the {\bf degree Rao bound}, and the sum polynomial
\(\Psi_e(x)=\sum_{\substack{u\in\mathcal{D}^*\\|u|\leq e(T)}}\Phi_u(x)\)
the {\bf degree Wilson polynomial}.

\begin{definition}
A $T$-design $Y$ in a multivariate $Q$-polynomial scheme will be said to be a \textbf{degree tight T-design} if it satisfies the degree Rao bound, i.e.
\[
|Y|=\sum_{\substack{\gamma\in\mathcal{D}^*\\|\gamma|\leq e(T)}}\mu_\gamma.
\]
\end{definition}

\subsection{Four important examples}
The following examples translate the abstract notion of a \(t\)-design in a multivariate \(Q\)-polynomial scheme into concrete combinatorial structures, laying the groundwork for further discussions.

\subsubsection{$q$-ary Johnson scheme $J_q(w,n)$}
In this example, we use the dual index weight \(d_{\mathcal D^*}(i,j)=i+j\).
By \cite[Lemma 5.1]{qary}, this weight is compatible with the Krein support. Hence the index set defining a distance-\(t\) design is \(T_t=\{(i,j)\in L:0<i+j\le t\}\). Assume that $t\leq m$. For each admissible $(r,s)\in L$ with $0<r+s\leq t$, set
\(T(r,s)=\{(i,j)\in L:i\leq r,\ j\leq s\}\setminus\{(0,0)\}\), then
\[
T_t=\bigcup_{\substack{(r,s)\in L\\0<r+s\leq t}}T(r,s).
\]
Indeed, if $(i,j)\in T_t$, then $(i,j)\in T(i,j)$; conversely, $(i,j)\in T(r,s)$ with $r+s\leq t$ implies $i+j\leq r+s\leq t$.

By \cite[Theorem~3.2]{qary}, since $r\leq r+s\leq t\leq m$, a $T(r,s)$-design is equivalent to an $(r,s)$-design in the sense of Nozaki-Watanabe. Consequently,
$Y\subseteq J_q(w,n)$ is a distance-$t$ design if and only if it is simultaneously an $(r,s)$-design for every admissible $(r,s)\in L$ satisfying $0<r+s\leq t$.

The combinatorial meaning is the following \cite[Definition~4.1]{qary}. For fixed $(r,s)$, choose $R\subseteq\{1,\ldots,n\}$ and $S\subseteq R$ with $|R|=r$ and $|S|=s$, and prescribe $\omega\in(\mathbb F_q^\ast)^S$. Define
\[
m_{R,S}(Y,\omega)=\bigl|\{y \in Y:R\subseteq\operatorname{supp}(y),\ y_i=\omega_i\\~\text{for all }i\in S\}\bigr|.
\]
Then the $(r,s)$-design condition means that $m_{R,S}(Y,\omega)$ depends only on $(r,s)$, and not on the particular choices of $R$, $S$, or $\omega$.

Thus, a distance-$t$ design in $J_q(w,n)$ is characterized combinatorially by the simultaneous uniformity of all support-incidence and nonzero-symbol conditions corresponding to \(0<r+s\leq t\).

\subsubsection{Sum-rank scheme}
Let \(\mathcal S=\bigotimes_{j=1}^{\ell}\Omega(n_j,m_j),~n_j\leq m_j\) be the product of bilinear forms schemes over $\mathbb F_q$, with vertex set
\(X=\prod_{j=1}^{\ell}\mathbb F_q^{\,n_j\times m_j}\).
For \(F=(F_1,\ldots,F_\ell),~G=(G_1,\ldots,G_\ell)\), the sum-rank distance is
\(d_{\rm SR}(F,G)=\sum_{j=1}^{\ell}\operatorname{rk}(F_j-G_j)\). Under the standard self-dual rank labeling of the bilinear forms schemes, the primitive idempotents are indexed by \(k=(k_1,\ldots,k_\ell),~0\leq k_j\leq n_j\), with
\(d_{\mathcal D^*}(k)=|k|=\sum_{j=1}^{\ell}k_j\). Hence \(T_t=\{k:1\leq |k|\leq t\}\), and a distance-$t$ design is precisely a $T_t$-design.

To describe its combinatorial meaning, identify $F_j\in\mathbb F_q^{n_j\times m_j}$ with the linear map \(L_{F_j}:\mathbb F_q^{n_j}\longrightarrow\mathbb F_q^{m_j},~L_{F_j}(x)=F_j^{T}x\). The lower semilattice of bilinear forms consists of pairs \((U_j,f_j),~U_j\leq\mathbb F_q^{n_j},~f_j:U_j\longrightarrow\mathbb F_q^{m_j}\)
with $f_j$ linear, ordered by \((U_j,f_j)\preceq(U'_j,f'_j)~\Longleftrightarrow~U_j\subseteq U'_j\ \text{and}\ f'_j|_{U_j}=f_j.\)
Its rank is $\dim U_j$ \cite[Section~7]{lattice}. Thus $(U_j,f_j)$ is incident with the
complete map $L_{F_j}$ exactly when \(L_{F_j}|_{U_j}=f_j\).

For the product semilattice, the height vector is \(u=(\dim U_1,\ldots,\dim U_\ell)\). By the product-poset characterization \cite[Theorem~2.3]{DP}, a subset $Y\subseteq X$ is a
distance-$t$ design if and only if, for every choice of $(U_j,f_j)$ satisfying \(\sum_{j=1}^{\ell}\dim U_j\leq t\), the number
\[
\left|\left\{(F_1,\ldots,F_\ell)\in Y:L_{F_j}|_{U_j}=f_j\text{ for all }j\right\}\right|
\]
depends only on the dimension vector \((\dim U_1,\ldots,\dim U_\ell)\), and not on the particular subspaces $U_j$ or linear maps $f_j$.

\subsubsection{Mixed alphabet scheme}
Consider the product Hamming scheme
\(\mathfrak{X}=H(n_1,q_1)\otimes H(n_2,q_2)\otimes\cdots\otimes H(n_m,q_m)\).
Under the standard self-dual labeling of the Hamming schemes, the primitive idempotents of \(\mathfrak{X}\) are indexed by \(j=(j_1,\ldots,j_m),~0\leq j_i\leq n_i\), and the corresponding dual distance is \(d_{\mathcal D^*}(j)=|j|=\sum_{i=1}^{m}j_i\). Hence
\(T_t=\{j:1\leq |j|\leq t\}\), and a distance-$t$ design is precisely a $T_t$-design.

For the $i$-th Hamming factor, a rank-$j_i$ element of the associated Hamming lattice is obtained by choosing $j_i$ coordinates and prescribing their symbols \cite[Example~5.2]{OOR}. Thus, in the product lattice, an element of height vector \(j=(j_1,\ldots,j_m)\) specifies $j_i$ coordinates in the $i$-th block, with a total of $|j|$ prescribed coordinates.

By the product-poset characterization \cite[Theorem~2.3]{DP}, a subset $Y\subseteq X$ is a
distance-$t$ design if and only if these prescribed coordinate patterns occur equally often for every height vector $j$ with $|j|\leq t$. Equivalently, when the elements of $Y$ are arranged as a \(|Y|\times(n_1+\cdots+n_m)\) array whose $i$-th block consists of $n_i$ columns over an alphabet of size $q_i$, then for every choice of $t$ columns, each tuple allowed by the alphabets of those columns occurs equally often among the rows. Thus, a distance-$t$ design is precisely a mixed-level orthogonal array of strength $t$.

\subsubsection{NRT scheme}
In the NRT scheme (the ordered Hamming scheme $\vec{H}(s,\ell,q)$) \cite{NR, NRT}, the primitive idempotents are indexed by shapes $e = (e_0, \dots, e_\ell)$ with
$\sum_{i=0}^\ell e_i = s$. Under the standard formally self-dual shape labeling, the corresponding dual distance is the height \(d_{\mathcal D^*}(e)=h(e)=\sum_{i=0}^{\ell} i e_i\). Hence \(T_t=\{\, e = (e_0, \dots, e_\ell) \mid \sum_{i=0}^\ell e_i = s,\;1 \le h(e) \le t \,\}\). Consequently, a subset $Y \subseteq X$ is a $t$-design if and only if it is a $T_t$-design.

By \cite[Theorem~3.4]{NRT}, the vanishing of the dual coefficients for all shapes of nonzero height at most $t$ is equivalent to the ordered orthogonal array condition, thus a $t$-design in the NRT scheme is equivalent to an ordered orthogonal array of strength $t$.

More explicitly, write the elements of $Y$ as the rows of a $|Y|\times s\ell$ array, with the columns partitioned into $s$ ordered blocks of size $\ell$. For any integers
\(0\leq t_j\leq \ell,~\sum_{j=1}^{s}t_j=t\), select the first $t_j$ columns from the $j$-th block. Then every $q$-ary $t$-tuple occurs equally often in the resulting $t$ columns.
\section{Conclusion and open problems}
In this paper we have strived to extend from Q-polynomial to multivariate Q-polynomial association schemes the main results of \cite[Chapter 5]{D2} on codes of given degree and designs of given strength.
Thus, to be specific, we have generalized Sections 5.3.1 and 5.3.2 of that thesis. A result that has resisted our efforts is \cite[Theorem 5.22]{D2} that requires multivariate analogues of Christoffel numbers. The theory of multivariate  Gaussian quadrature is very involved and difficult to apply \cite{X}. This is the first problem open by this work.  An intriguing  question is the combinatorial interpretation of designs in the association schemes corresponding to the homogeneous and Lee distance \cite{MW,S3}.
\section*{Appendix: Self-dual association schemes}
\subsection*{Direct Product of Association Schemes}
Following the approach in \cite{DP} we describe direct products.
For $1 \leq i \leq m$, let $(Y_i, \mathcal{A}_i)$ be a $d_i$-class association scheme with adjacency matrices $\mathcal{A}_i$. The direct product of these schemes \cite{DP} is the association scheme
\[
(X, \mathcal{A}) = (Y_1, \mathcal{A}_1) \otimes (Y_2, \mathcal{A}_2) \otimes \cdots \otimes (Y_m, \mathcal{A}_m)
\]
defined by
\[
X = Y_1 \times Y_2 \times \cdots \times Y_m
\]
and
\[
\mathcal{A} = \left\{ \bigotimes_{i=1}^m M_i : M_i \in \mathcal{A}_i, 1 \leq i \leq m \right\},
\]
where
\[
\bigotimes_{i=1}^m M_i = M_1 \otimes M_2 \otimes \cdots \otimes M_m
\]
is the $m$-fold Kronecker product of matrices. The fact that this always gives an association scheme follows from the properties
\[
\left( \bigotimes_{i=1}^m M_i \right) \left( \bigotimes_{i=1}^m N_i \right) = \bigotimes_{i=1}^m (M_i N_i),~
\left( \bigotimes_{i=1}^m M_i \right) \circ \left( \bigotimes_{i=1}^m N_i \right) = \bigotimes_{i=1}^m (M_i \circ N_i)
\]
of the Kronecker product. From the first identity, it follows that the $P$-matrix for this product scheme is given by the Kronecker product of the $P$-matrices for the component schemes $(Y_i, \mathcal{A}_i)$. Similarly, the second change-of-basis matrix $Q$ for the product scheme is given by $\bigotimes_{i=1}^m Q_i$ where $Q_i$ is the $Q$-matrix for the $i^{\text{th}}$ component scheme.

For the purpose of the following proposition, each self-dual association scheme is assumed to be equipped with a compatible labeling of its relations and primitive idempotents under which its $P$-matrix equals its $Q$-matrix.

{ \prop\label{dp}
The direct product of finitely many self-dual association schemes is also self-dual.
}

\begin{proof}
Let $(X, \mathcal{A}) = \bigotimes_{i=1}^m (Y_i, \mathcal{A}_i)$ be the direct product of self-dual schemes $(Y_i, \mathcal{A}_i)$. By the properties recalled above, the $P$-matrix and $Q$-matrix of the product scheme are
\[
P = \bigotimes_{i=1}^m P_i \quad \text{and} \quad Q = \bigotimes_{i=1}^m Q_i.
\]
Since each $(Y_i, \mathcal{A}_i)$ is self-dual, $P_i = Q_i$ for all $i$. Substitution yields $P = Q$, hence $(X, \mathcal{A})$ is self-dual.
\end{proof}

\subsection*{Characterization of self-dual translation schemes}
{ Let $X$ be a finite abelian group. Assume $X$ is the groundset of a translation association scheme with $s$ classes, then there exists a {\bf partition} $P_0=\{0\},P_1,\dots,P_s$ of $X$ such that the relations of the scheme are given by $x R_i y$ iff $x-y \in P_i$. We are  going to characterize the partitions that turn $(X,R)$ into a self-dual association scheme. Denote the indicator function of \(P_i \in P\) by \(\chi_i\):
\[
\chi_i(x) =
\begin{cases}
	1 & \text{if } x \in P_i, \\
	0 & \text{if } x \notin P_i.
\end{cases}
\]
Let \( \mathbb{C}^X \) denote the family of all complex functions \( X \to \mathbb{C} \). For any given partition \( P \), we define
\[
L(P) \triangleq \{ f \in \mathbb{C}^X : f(x) = \sum_{i=0}^s f_i \chi_i(x) \},
\]
where \( f_i \in \mathbb{C}, i = 0, 1, \dots, s \).
\begin{definition}
Let \(X\) be a finite abelian group and \(P = \{P_0 , P_1 , \dots , P_s\}\) be a partition of \(X\) with \(P_0 = \{0\}\).  Thus $X=\bigsqcup_{i=0}^sP_i.$
Let \(L^*(P) = \{ f^* : f \in L(P) \}\) be the image of \(L(P)\) under the Fourier transform \(f^*(y) = \sum_{x \in X} \Psi_x(y) f(x)\), where \(\Psi_x\) are the additive characters of \(X\).
If \(L^*(P) = L(P)\), then \(P\) is called an \textbf{F-partition}.
\end{definition}
\cite[Theorem 2]{ZE} reveals that the essence of the self-duality of a translation association scheme is the invariance of the partition $P$ under the Fourier transform:
$(X,R)$ is self-dual iff $P$ is an $F$-partition.
On the other hand,  \cite[Corollary 7]{TAS} tells us when the translation association scheme generated by the $G$-orbits partition ( for some group $G$ acting on $X$)  is self-dual. The existence of an adjoint map actually guarantees that the orbit partition satisfies Fourier invariance, and the orbit partition in \cite{TAS} is a special type of $F$-partition in the sense of \cite{ZE}.
}

{
The reference \cite{MW} constructs another class of $F$-partitions that are not generated by group orbits (e.g., the four-block partition for the homogeneous weight) and verifies that they also satisfy the Fourier invariance condition of \cite{ZE}. Therefore, \cite{MW} and \cite{TAS} provide two different sources of $F$-partitions: the former comes from coarsest constructions for concrete metrics, whereas the latter comes from group actions with an adjoint map.
}
\subsection*{Examples}
We now illustrate five well-known families of translation schemes, all of which are self-dual.

\begin{enumerate}
\item \textbf{Hamming scheme \(\mathrm{H}(n,q)\).}
The Hamming scheme \(\mathrm{H}(n,q)\) enjoys the MacWilliams relations on the complete weight enumerator.
It is well-known to be self-dual \cite{D}.

\item \textbf{Mixed alphabet scheme.}
A mixed alphabet scheme is essentially a direct product of Hamming schemes \(\mathrm{H}(n,q)\) over possibly different $n$'s and $q$'s. As each factor is self-dual, the product remains self-dual
by Proposition \ref{dp}.
	
\item \textbf{Lee scheme.}
By \cite{D2} the Lee scheme is self-dual.
	
\item \textbf{Homogeneous metric scheme.}
The homogeneous metric scheme is self-dual under the natural character identification. Indeed, the character-sum computation for the homogeneous weight-unit partition in \cite{MW} shows that its dual partition has the same blocks \(Z,U,S,R\).
	
\item \textbf{Sum-rank scheme.}
The sum-rank scheme can be viewed as a direct product of bilinear forms schemes. The bilinear forms scheme is self-dual by \cite[(3.9)]{D3}, and self-duality is preserved under direct products, consequently the sum-rank scheme is self-dual by Proposition \ref{dp}.
\end{enumerate}

\section*{Acknowledgement:} The authors are grateful to Eiichi Bannai for helpful discussions.

\section*{Declaration of AI use}
During the preparation and revision of this article, the authors used ChatGPT (OpenAI) to assist with language refinement, the exploration and development of ideas, and improvements to the clarity. The authors reviewed and revised the AI-assisted content and take full responsibility for the results, the arguments, and the final form of the article.

\end{document}